\documentclass[11pt]{amsart}
\usepackage{amsmath}
\usepackage{amsxtra}
\usepackage{amscd}
\usepackage{amsthm}
\usepackage{amsfonts}
\usepackage{amssymb}
\usepackage{eucal}

\newtheorem{lem}[subsection]{Lemma}
\newtheorem{prop}[subsection]{Proposition}

\newtheorem{thm}[subsection]{Theorem}

\newtheorem{rem}[subsection]{Remark}

\theoremstyle{definition}

\theoremstyle{remark}

\newcommand{\nc}{\newcommand}
\nc{\renc}{\renewcommand}
\nc{\ssec}{\subsection}
\nc{\sssec}{\subsubsection}
\nc{\on}{\operatorname}

\nc\ol{\overline}
\nc\ul{\underline}
\nc\wt{\widetilde}
\nc\tboxtimes{\wt{\boxtimes}}
\nc{\alp}{\alpha}

\nc{\ZZ}{{\mathbb Z}}
\nc{\NN}{{\mathbb N}}
\nc{\CC}{{\mathbb C}}
\nc{\OO}{{\mathbb O}}
\renc{\SS}{{\mathbb S}}
\nc{\DD}{{\mathbb D}}
\nc{\GG}{{\mathbb G}}

\nc{\Fq}{{\mathbb F}_q}
\nc{\Fqb}{\ol{{\mathbb F}_q}}
\nc{\Ql}{\ol{{\mathbb Q}_\ell}}
\nc{\id}{\text{id}}
\nc\X{\mathcal X}

\nc{\Hom}{\on{Hom}}
\nc{\Lie}{\on{Lie}}
\nc{\Loc}{\on{Loc}}
\nc{\Pic}{\on{Pic}}
\nc{\Bun}{\on{Bun}}
\nc{\IC}{\on{IC}}
\nc{\Aut}{\on{Aut}}
\nc{\rk}{\on{rk}}
\nc{\Sh}{\on{Sh}}
\nc{\Perv}{\on{Perv}}
\nc{\pos}{{\on{pos}}}
\nc{\Conv}{\on{Conv}}
\nc{\Sph}{\on{Sph}}
\nc{\Sym}{\on{Sym}}
\nc{\BunBb}{\overline{\Bun}_B}
\nc{\Buno}{\overset{o}{\Bun}}
\nc{\BunPb}{{\overline{\Bun}_P}}
\nc{\BunBM}{\overline{\Bun}_{B(M)}}
\nc{\BunPbw}{{\widetilde{\Bun}_P}}
\nc{\BunBP}{\widetilde{\Bun}_{B,P}}
\nc{\GUb}{\overline{G/U}}
\nc{\GUPb}{\overline{G/U(P)}}

\nc{\Hhom}{\underline{\on{Hom}}}
\nc\syminfty{\on{Sym}^{\infty}}
\nc\lal{\ol{\lambda}}
\nc\xl{\ol{x}}
\nc\thl{\ol{\theta}}
\nc\nul{\ol{\nu}}
\nc\mul{\ol{\mu}}
\nc\Sum\Sigma
\nc{\oX}{\overset{o}{X}{}}

\nc{\M}{{\mathcal M}}
\nc{\N}{{\mathcal N}}
\nc{\F}{{\mathcal F}}
\nc{\D}{{\mathcal D}}
\nc{\Q}{{\mathcal Q}}
\nc{\Y}{{\mathcal Y}}
\nc{\G}{{\mathcal G}}
\nc{\E}{{\mathcal E}}
\nc{\CalC}{{\mathcal C}}
\nc\Dh{\widehat{\D}}
\renewcommand{\O}{{\mathcal O}}
\nc{\C}{{\mathcal C}}
\nc{\K}{{\mathcal K}}
\renewcommand{\H}{{\mathcal H}}

\nc{\T}{{\mathcal T}}
\nc{\V}{{\mathcal V}}
\renc{\P}{{\mathcal P}}
\nc{\A}{{\mathcal A}}
\nc{\B}{{\mathcal B}}
\nc{\U}{{\mathcal U}}

\nc{\Gr}{\on{Gr}}

\nc{\frn}{{\check{\mathfrak u}(P)}}
\nc\f{{\mathfrak f}}

\nc{\q}{{\mathfrak q}}
\nc{\p}{{\mathfrak p}}
\nc{\s}{{\mathfrak s}}
\nc\w{\text{w}}

\nc\mathi\iota
\nc\Spec{\on{Spec}}
\nc\Mod{\on{Mod}}
\nc{\tw}{\widetilde{\mathfrak t}}
\nc{\pw}{\widetilde{\mathfrak p}}
\nc{\qw}{\widetilde{\mathfrak q}}
\nc{\jw}{\widetilde j}

\nc{\grb}{\overline{\Gr}}
\nc{\I}{\mathcal I}

\nc{\lambdach}{{\check\lambda}}
\nc{\Lambdach}{{\check\Lambda}{}}
\nc{\much}{{\check\mu}}
\nc{\omegach}{{\check\omega}}
\nc{\nuch}{{\check\nu}}
\nc{\etach}{{\check\eta}}
\nc{\alphach}{{\check\alpha}}
\nc{\betach}{{\check\beta}}
\nc{\rhoch}{{\check\rho}}
\nc{\ch}{{\check h}}

\nc{\Hb}{\overline{\H}}

\nc{\BA}{{\mathbb{A}}}
\nc{\BC}{{\mathbb{C}}}
\nc{\BM}{{\mathbb{M}}}
\nc{\BN}{{\mathbb{N}}}
\nc{\BP}{{\mathbb{P}}}
\nc{\BR}{{\mathbb{R}}}
\nc{\BZ}{{\mathbb{Z}}}
\nc{\BS}{{\mathbb{S}}}

\nc{\CA}{{\mathcal{A}}}
\nc{\CB}{{\mathcal{B}}}
\nc{\CE}{{\mathcal{E}}}
\nc{\CF}{{\mathcal{F}}}
\nc{\CG}{{\mathcal{G}}}
\nc{\CH}{{\mathcal{H}}}
\nc{\CL}{{\mathcal{L}}}
\nc{\CM}{{\mathcal{M}}}
\nc{\CN}{{\mathcal{N}}}
\nc{\CO}{{\mathcal{O}}}
\nc{\CP}{{\mathcal{P}}}
\nc{\CQ}{{\mathcal{Q}}}
\nc{\CR}{{\mathcal{R}}}
\nc{\CS}{{\mathcal{S}}}
\nc{\CT}{{\mathcal{T}}}
\nc{\CU}{{\mathcal{U}}}
\nc{\CV}{{\mathcal{V}}}
\nc{\CW}{{\mathcal{W}}}
\nc{\CZ}{{\mathcal{Z}}}

\nc{\cM}{{\check{\mathcal M}}{}}
\nc{\csM}{{\check{\mathcal A}}{}}
\nc{\oM}{{\overset{\circ}{\mathcal M}}{}}
\nc{\obM}{{\overset{\circ}{\mathbf M}}{}}
\nc{\oCA}{{\overset{\circ}{\mathcal A}}{}}
\nc{\obA}{{\overset{\circ}{\mathbf A}}{}}
\nc{\ooM}{{\overset{\circ}{M}}{}}
\nc{\osM}{{\overset{\circ}{\mathsf M}}{}}
\nc{\vM}{{\overset{\bullet}{\mathcal M}}{}}
\nc{\nM}{{\underset{\bullet}{\mathcal M}}{}}
\nc{\oD}{{\overset{\circ}{\mathcal D}}{}}
\nc{\obD}{{\overset{\circ}{\mathbf D}}{}}
\nc{\oA}{{\overset{\circ}{\mathbb A}}{}}
\nc{\op}{{\overset{\bullet}{\mathbf p}}{}}
\nc{\cp}{{\overset{\circ}{\mathbf p}}{}}
\nc{\oU}{{\overset{\bullet}{\mathcal U}}{}}
\nc{\oZ}{{\overset{\circ}{\mathcal Z}}{}}
\nc{\ofZ}{{\overset{\circ}{\mathfrak Z}}{}}

\nc{\fa}{{\mathfrak{a}}}
\nc{\fb}{{\mathfrak{b}}}
\nc{\fg}{{\mathfrak{g}}}
\nc{\fgl}{{\mathfrak{gl}}}
\nc{\fh}{{\mathfrak{h}}}
\nc{\fri}{{\mathfrak{i}}}
\nc{\fj}{{\mathfrak{j}}}
\nc{\fm}{{\mathfrak{m}}}
\nc{\fn}{{\mathfrak{n}}}
\nc{\ft}{{\mathfrak{t}}}
\nc{\fu}{{\mathfrak{u}}}
\nc{\fz}{{\mathfrak{z}}}
\nc{\fp}{{\mathfrak{p}}}
\nc{\frr}{{\mathfrak{r}}}
\nc{\fs}{{\mathfrak{s}}}
\nc{\fsl}{{\mathfrak{sl}}}
\nc{\hsl}{{\widehat{\mathfrak{sl}}}}
\nc{\hgl}{{\widehat{\mathfrak{gl}}}}
\nc{\hg}{{\widehat{\mathfrak{g}}}}
\nc{\chg}{{\widehat{\mathfrak{g}}}{}^\vee}
\nc{\hn}{{\widehat{\mathfrak{n}}}}
\nc{\chn}{{\widehat{\mathfrak{n}}}{}^\vee}

\nc{\fA}{{\mathfrak{A}}}
\nc{\fB}{{\mathfrak{B}}}
\nc{\fD}{{\mathfrak{D}}}
\nc{\fE}{{\mathfrak{E}}}
\nc{\fF}{{\mathfrak{F}}}
\nc{\fG}{{\mathfrak{G}}}
\nc{\fK}{{\mathfrak{K}}}
\nc{\fL}{{\mathfrak{L}}}
\nc{\fM}{{\mathfrak{M}}}
\nc{\fN}{{\mathfrak{N}}}
\nc{\frP}{{\mathfrak{P}}}
\nc{\fU}{{\mathfrak{U}}}
\nc{\fZ}{{\mathfrak{Z}}}

\nc{\bb}{{\mathbf{b}}}
\nc{\bc}{{\mathbf{c}}}
\nc{\be}{{\mathbf{e}}}
\nc{\bj}{{\mathbf{j}}}
\nc{\bn}{{\mathbf{n}}}
\nc{\bp}{{\mathbf{p}}}
\nc{\bq}{{\mathbf{q}}}
\nc{\bfu}{{\mathbf{u}}}
\nc{\bv}{{\mathbf{v}}}
\nc{\bx}{{\mathbf{x}}}
\nc{\by}{{\mathbf{y}}}
\nc{\bw}{{\mathbf{w}}}
\nc{\bA}{{\mathbf{A}}}
\nc{\bB}{{\mathbf{B}}}
\nc{\bC}{{\mathbf{C}}}
\nc{\bD}{{\mathbf{D}}}
\nc{\bF}{{\mathbf{F}}}
\nc{\bH}{{\mathbf{H}}}
\nc{\bK}{{\mathbf{K}}}
\nc{\bM}{{\mathbf{M}}}
\nc{\bN}{{\mathbf{N}}}
\nc{\bO}{{\mathbf{O}}}
\nc{\bT}{{\mathbf{T}}}
\nc{\bV}{{\mathbf{V}}}
\nc{\bW}{{\mathbf{W}}}
\nc{\bX}{{\mathbf{X}}}
\nc{\bP}{{\mathbf{P}}}
\nc{\bZ}{{\mathbf{Z}}}

\nc{\sA}{{\mathsf{A}}}
\nc{\sB}{{\mathsf{B}}}
\nc{\sC}{{\mathsf{C}}}
\nc{\sD}{{\mathsf{D}}}
\nc{\sF}{{\mathsf{F}}}
\nc{\sK}{{\mathsf{K}}}
\nc{\sM}{{\mathsf{M}}}
\nc{\sO}{{\mathsf{O}}}
\nc{\sQ}{{\mathsf{Q}}}
\nc{\sP}{{\mathsf{P}}}
\nc{\sT}{{\mathsf{T}}}
\nc{\sZ}{{\mathsf{Z}}}
\nc{\sfp}{{\mathsf{p}}}
\nc{\sr}{{\mathsf{r}}}
\nc{\sfb}{{\mathsf{b}}}
\nc{\sfc}{{\mathsf{c}}}
\nc{\sd}{{\mathsf{d}}}

\nc{\BK}{{\bar{K}}}

\nc{\tA}{{\widetilde{\mathbf{A}}}}
\nc{\tB}{{\widetilde{\mathcal{B}}}}
\nc{\tg}{{\widetilde{\mathfrak{g}}}}
\nc{\tG}{{\widetilde{G}}}
\nc{\TM}{{\widetilde{\mathbb{M}}}{}}
\nc{\tO}{{\widetilde{\mathsf{O}}}{}}
\nc{\tU}{{\widetilde{\mathfrak{U}}}{}}
\nc{\TZ}{{\tilde{Z}}}
\nc{\tx}{{\tilde{x}}}
\nc{\tbv}{{\tilde{\bv}}}
\nc{\tfP}{{\widetilde{\mathfrak{P}}}{}}
\nc{\tz}{{\tilde{\zeta}}}
\nc{\tmu}{{\tilde{\mu}}}

\nc{\urho}{\underline{\rho}}
\nc{\uB}{\underline{B}}
\nc{\uC}{{\underline{\mathbb{C}}}}
\nc{\ui}{\underline{i}}
\nc{\uj}{\underline{j}}
\nc{\ofP}{{\overline{\mathfrak{P}}}}
\nc{\oB}{{\overline{\mathcal{B}}}}
\nc{\og}{{\overline{\mathfrak{g}}}}
\nc{\oI}{{\overline{I}}}

\nc{\eps}{\varepsilon}
\nc{\hrho}{{\hat{\rho}}}

\nc{\one}{{\mathbf{1}}}
\nc{\two}{{\mathbf{t}}}

\nc{\Rep}{{\mathop{\operatorname{\rm Rep}}}}
\nc{\Tot}{{\mathop{\operatorname{\rm Tot}}}}
\nc{\Ker}{{\mathop{\operatorname{\rm Ker}}}}
\nc{\Hilb}{{\mathop{\operatorname{\rm Hilb}}}}
\nc{\End}{{\mathop{\operatorname{\rm End}}}}
\nc{\Ext}{{\mathop{\operatorname{\rm Ext}}}}
\nc{\CHom}{{\mathop{\operatorname{{\mathcal{H}}\it om}}}}
\nc{\GL}{{\mathop{\operatorname{\rm GL}}}}
\nc{\gr}{{\mathop{\operatorname{\rm gr}}}}
\nc{\Id}{{\mathop{\operatorname{\rm Id}}}}
\nc{\defi}{{\mathop{\operatorname{\rm def}}}}
\nc{\length}{{\mathop{\operatorname{\rm length}}}}
\nc{\supp}{{\mathop{\operatorname{\rm supp}}}}

\nc{\Cliff}{{\mathsf{Cliff}}}
\nc{\Fl}{{\mathsf{Fl}}}
\nc{\Fib}{{\mathsf{Fib}}}
\nc{\Coh}{{\mathsf{Coh}}}
\nc{\FCoh}{{\mathsf{FCoh}}}

\nc{\reg}{{\text{\rm reg}}}

\nc{\cplus}{{\mathbf{C}_+}}
\nc{\cminus}{{\mathbf{C}_-}}
\nc{\cthree}{{\mathbf{C}_*}}
\nc{\Qbar}{{\bar{Q}}}

\nc{\bh}{{\bar{h}}}
\nc{\bOmega}{{\overline{\Omega}}}

\nc{\seq}[1]{\stackrel{#1}{\sim}}

\nc{\BBB}{{\mathbf B}}
\nc{\BBC}{{\mathbf C}}
\nc{\BBD}{{\mathbf D}}
\nc{\BBK}{{\mathbf K}}

\nc{\fV}{{\mathfrak{V}}}
\nc{\BG}{{\mathbb{G}}}
\nc{\oF}{{\overset{\circ}{\fF}}}

\theoremstyle{remark}
\newtheorem{Rem}{Remark}

\newcommand{\cal}{\mathcal}

\newcommand{\bu}{\bullet}

\newcommand{\To}{\longrightarrow}
\newcommand{\iso}{{\widetilde \longrightarrow}}

\renewcommand{\P}{{\cal P}}

\nc{\fBt}{{\overset{\bullet}{\mathfrak{B}}}{}}
\nc{\fBtil}{\tilde{\mathfrak{B}}}
\nc{\varpit}{{\overset{\bullet}{\varpi}}{}}

\renewcommand{\H}{{\cal H}}

\newcommand{\oplusl}{\bigoplus\limits}
\newcommand{\cupl}{\bigcup\limits}

\newcommand{\g}{{\frak g}}

\newcommand{\Lg}{{\check{\mathfrak g}}{}}

\newcommand{\Lt}{{\check{\mathfrak t}}{}}

\newcommand{\LG}{{\check G}{}}

\newcommand{\LR}{{\check R}{}}

\newcommand{\LT}{{\check T}{}}

\renewcommand{\O}{{\cal O}}

\newcommand{\GO}{{G({\mathbf O})}}
\newcommand{\TO}{{T({\mathbf O})}}
\newcommand{\LGO}{{{\check G}({\mathbf O})}}
\newcommand{\GK}{{G({\mathbf F})}}
\newcommand{\GF}{{G({\mathbf F})}}
\newcommand{\TF}{{T({\mathbf F})}}

\newcommand{\Gm}{{{\mathbb G}_m}}

\newcommand{\beq}{\begin{equation}}
\newcommand{\eeq}{\end{equation}}

\renewcommand{\proof}{{\it Proof }}

\begin{document}

\author{Roman Bezrukavnikov, Michael Finkelberg, and Ivan Mirkovi\'c}
\title{Equivariant ($K$-)homology of affine Grassmannian and Toda lattice}
\address{{\it Address}:\newline
R.B.: Dept. of Math., Massachusetts Institute of Technology, 77 Massachusetts ave, Cambridge MA,
02139 USA; \newline
M.F.: Independent Moscow Univ., 11 Bolshoj Vlasjevskij per., 
Moscow 119002, Russia; \newline
I.M.: Dept. of Math., The University of Massachusetts, Amherst, MA 01003, USA}

\email{\newline
bezrukav@math.mit.edu; fnklberg@mccme.ru;
mirkovic@math.umass.edu}

\maketitle


\section{Introduction}

\subsection{}
\label{Ginzburg}
Let $G$ be an almost simple complex algebraic group, and let $\Gr_G$ be
its affine Grassmannian. Recall that if we set
$\bO=\BC[[{\mathsf t}]],\ \bF=\BC(({\mathsf t}))$,
then $\Gr_G=\GK/\GO$.

 It is well-known that the subgroup $\Omega K$
of polynomial loops into a maximal compact subgroup $K\subset G$ 
projects isomorphically to $\Gr_G$; thus $\Gr_G$ acquires the structure 
of a topological group. An algebro-geometric counterpart of this structure
is provided by {\it the convolution diagram} $\GK\times_{\GO} \Gr_G\to \Gr_G$.


It allows one to define the {\it convolution} of two $\GO$ equivariant
 geometric objects (such as sheaves, or constrictible functions) on
$\Gr_G$. A famous example of such a structure
is the category of $\GO$ equivariant perverse sheaves on $\Gr$ (``Satake
category'' in the terminology of Beilinson and Drinfeld); this is a semi-simple
abelian category, and convolution provides it with a symmetric
monoidal structure. By results of \cite{g}, \cite{mv}, \cite{bd}
this category is identified with the category of
(algebraic) representations of the Langlands dual group.

The starting point for the present work was the observation that a similar
definition works in another setting, yielding a monoidal structure
 on the category
of $\GO$ equivariant {\it perverse coherent sheaves} on $\Gr$
(the ``coherent Satake category''). The latter is
a non-semisimple artinian abelian category, the heart of the middle
perversity  $t$-structure on the derived category of $\GO$ equivariant
coherent sheaves on $\Gr_G$; existence of this $t$-structure is due to the 
fact that dimensions of all $\GO$-orbits inside a given component 
of $\Gr_G$ are of the same parity, cf. \cite{b}.  
The resulting monoidal  category turns out to be non-symmetric, though 
its Grothendieck ring $K^\GO(\Gr_G)$ is commutative. One of the results
of this paper is a computation of this ring.
 Along with $K^\GO(\Gr_G)$ we compute its
 ``graded version'',
the ring $H^\GO(\Gr)$ of equivariant homology of $\Gr$, where
the algebra structure
is again provided by convolution.\footnote{The two rings are related via
the Chern character homomorphism from $K^\GO(\Gr)$ to the
completion of $H^{\GO}(\Gr)$.} 
(The ring $H^\GO_\bullet(\Gr_G)$ was essentially computed by
Dale Peterson ~\cite{p}, cf. also ~\cite{k2}.)

To describe the answer, let $\LG$ be the Langlands dual group to $G$,
and let $\Lg$ be its Lie algebra.
Consider the {\em universal centralizers} $\fZ^{\check G}_\Lg$
and $\fZ^{\check G}_\LG$: if we denote by $C_{\LG,\Lg}\subset\LG\times\Lg$
(resp. $C_{\LG,\LG}\subset\LG\times\LG$) the locally closed subvariety
formed by all the pairs $(g,x)$ such that $Ad_g(x)=x$ and $x$ is regular
(resp. all the pairs $(g_1,g_2)$ such that $Ad_{g_1}g_2=g_2$ and $g_2$
is regular), then $\fZ^{\check G}_\Lg$ (resp. $\fZ^{\check G}_\LG$) is the
categorical quotient $C_{\LG,\Lg}//\LG$ (resp. $C_{\LG,\LG}//\LG$) with
respect to the diagonal adjoint action of $\LG$.

We identify $\on{Spec}\left(H^\GO_\bullet(\Gr_G)\right)$ with
$\fZ^{\check G}_\Lg$.
Also, we identify $\on{Spec}\left(K^\GO(\Gr_G)\right)$ with {\em a variant} of
$\fZ^{\check G}_\LG$ (the isomorphism
$\on{Spec}\left(K^\GO(\Gr_G)\right)\simeq\fZ^{\check G}_\LG$ holds true iff
$G$ is of type $E_8$).

Notice that $\fZ^{\check G}_\Lg$ inherits a canonical symplectic structure
as a hamiltonian reduction of the cotangent bundle $\sT^*\LG$.
Also, $\fZ^{\check G}_\LG$ inherits a canonical Poisson structure as a
q-Hamiltonian reduction of the q-Hamiltonian $\LG$-space {\em internal
fusion double} $\bD(\LG)$ (see ~\cite{amm}); this Poisson structure is in
fact symplectic iff $\LG$ is simply connected (that is, $G$ is adjoint).

The corresponding 
 Poisson structures on $K^{\GO}(\Gr_G)$, $H^\GO(\Gr_G)$ come from
a deformation of these commutative algebras to non-commutative algebras
 $H^{\GO\ltimes\Gm}_\bullet(\Gr_G)$
(resp. $K^{\GO\ltimes\Gm}(\Gr_G)$); here $\Gm$ acts on $\Gr_G$
by loop rotation. We conjecture that the non-commutative algebra
 $H^{\GO\ltimes\Gm}_\bullet(\Gr_G)$ can also be obtained from
the ring of differential operators on $\LG$ by quantum Hamiltonian reduction.

The space $\fZ^{\check G}_\Lg$ contains an open piece $\fZ(\LG)$
which for $\LG$ adjoint (that is, for
$G$ simply connected) is a complexification of the Kostant's
phase space of the classical Toda lattice ~(\cite{k}, ~Theorem ~2.6).
We remark in passing that
Toda lattice also appears in the (apparently related) computations by
Givental, Kim and others of quantum 
cohomology of flag varietites (see e.g. \cite{kim}).

Our computation should be compared with (and is to a large extent inspired by)
\cite{g} where equivariant cohomology $H_\GO(\Gr_G)$ were 
computed\footnote{Another description for $H_\GO(\Gr_G)$ is provided
by a general result of \cite{kuku1}; in fact, its extension
from \cite{kuku2} gives also an answer for $K^\GO(\Gr_G)$,
and a similar technique can be applied to compute $H^\GO(\Gr_G)$.
However, this form of the answer does not make the relation to the (dual)
group geometry explicit.}
in terms of the $\LG$. (The precise relation between the two computations
is spelled out in Remark \ref{relation_with_ginzburg}).

\medskip

The second main object considered in the paper is another derived
category of coherent sheaves with a convolution monoidal structure,
namely the derived category  $D^bCoh^\GO_{\Lambda_G}(T^*\Gr)$
of
$\GO$-equivariant coherent sheaves on the cotangent bundle of $\Gr_G$
supported on the union $\Lambda_G$ of conormal bundles to the $\GO$-orbits
(the definition of involved objects requires extra work since $\Gr_G$
is infinite dimensional).
(In this case we do not find a $t$-structure compatible with convolution,
so all we get is a monoidal triangulated category).
Notice that the singular support of a $\GO$-equivariant $D$-module on $\Gr_G$
is an object of $Coh^\GO_{\Lambda_G}(T^*\Gr)$, thus this category
can be considered a ``classical limit'' of the (derived) Satake category.
We compute the Grothendieck ring of  $D^bCoh^\GO_{\Lambda_G}(T^*\Gr)$
 identifying its spectrum with $(T\times \LT)/W$,
where $T\subset G$, and $\LT\subset \LG$ are Cartan subgroups.
This is a singular variety birationally equivalent to 
 $\on{Spec}\left(K^\GO(\Gr_G)\right)$. Unlike the latter, the former
remains unchanged if we replace $G$ by $\LG$. This motivates
a conjecture that the corresponding triangulated monoidal categories
for $G$ and $\LG$ are equivalent. The conjecture is compatible with 
a ``classical limit'' of the geometric Langlands conjecture of
Beilinson and Drinfeld (see \ref{classical_langlands}
 below for a more precise statement of the conjecture).

Finally, we remark that  the convolution of
$\GO$-equivariant perverse coherent sheaves is closely related to the
{\em fusion product} of $\GO$-modules introduced by
B.~Feigin\footnote{The relation between convolution and fusion was known
to B.~Feigin since 1997.} ~\cite{fl}
(see Section ~\ref{seven}).
In fact, our desire to understand the category $\CP^\GO(\Gr_G)$, and
the work ~\cite{fl} of B.~Feigin and S.~Loktev, was one of the motivations for
the present work.

\subsection{Acknowledgments}
We are obliged to D.~Gaitsgory, V.~Ginzburg,
D.~Kazhdan and S.~Loktev for help with the references, and especially to
P.~Etingof for the explanations about integrable systems.
R.B. is grateful to the Independent Moscow University and the 
American Embassy in Moscow
for granting him an opportunity to complete the work on this paper;
he is partially supported by NSF grant  DMS-0071967.
M.F. is grateful to the University of Massachusetts
at Amherst and the University of Chicago for the hospitality and support. 
His research was conducted for the Clay Mathematical Institute and partially
supported by the CRDF award RM1-2545-MO-03.

\bigskip

{\em This  update corrects a mistake appearing in the version of the article published in 2005, and a mathematical typo
in the previous arxiv version.} The modifications
include key definitions in sections~\ref{KosSte}--\ref{Aff} and the proof of Proposition~\ref{razdutija} in section~\ref{proof28}.
We are very grateful to Yunsong Wei for pointing out a mistake and to C\'edric
Bonnaf\'e and Eugene Gorsky for  helpful discussions. We are also thank Tom Gannon
for pointing out a discrepancy in the previous arxiv version.

\section{Notations and statements of the results}

\subsection{Kostant slices and Steinberg sections}\label{KosSte}
$G$ is an almost simple algebraic group with the Lie algebra
$\fg$. We choose a principal $\fsl_2$ triple $(e,h,f)$ in $\fg$. Let
$\phi\colon \fsl_2\to\fg$
be the corresponding homomorphism.
We denote by $\fz(e)$ the centralizer of $e$ in $\fg$.
We denote by $\Sigma_\fg\subset\fg$
the {\em Kostant slice} $\fz(e)+f$.
It is known that $\Sigma_\fg\subset\fg^{reg}$
and the projection to the categorical
quotient $\Sigma_\fg\hookrightarrow\fg\twoheadrightarrow\fg//Ad_G=\ft/W$
induces an isomorphism $\Sigma_\fg\simeq\ft/W$.

If $G$ is simply connected, we choose a Cartan torus $T$, and a Borel subgroup $B$ containing $T$,
with the unipotent radical $U$. The opposite Borel subgroup and its unipotent radical are denoted
by $B_-,U_-$ respectively. We also choose a representative $c$ in the normalizer $N_G(T)$
of a Coxeter element in the Weyl group $W(G,T)$. The {\em Steinberg cross-section} $\Sigma_G$
is defined as $U_-c\cap cU$. It is known that $\Sigma_G\subset G^{reg}$, and the composed morphism
$\Sigma_G\hookrightarrow G\to G//Ad_G$ is an isomorphism. The Steinberg cross-section is
isomorphic to an affine space. If $G$ is not necessarily simply connected, we define $\Sigma_G$
as the image in $G$ of the Steinberg cross-section of its universal cover $G^{sc}$. In general,
$\Sigma_G$ is not necessarily smooth.

\subsection{The universal centralizers}
We consider the locally closed subvariety $C_{\fg,\fg}\subset\fg\times\fg$
(resp. $C_{\fg,G}\subset\fg\times G,\ C_{G,\fg}\subset G\times\fg,\
C_{G,G}\subset G\times G$) formed by all the pairs $(x_1,x_2)$ such that
$[x_1,x_2]=0$ and $x_2$ is regular
(resp. all the pairs $(x,g)$ such that $Ad_g(x)=x$ and $g$ is regular;
all the pairs $(g,x)$ such that $Ad_g(x)=x$ and $x$ is regular;
all the pairs $(g_1,g_2)$ such that $Ad_{g_1}(g_2)=g_2$ and $g_2$ is regular).
The categorical quotients with respect to the diagonal adjoint action of $G$
are denoted respectively $C_{\fg,\fg}//G=\fZ_\fg^\fg,\ C_{\fg,G}//G=\fZ_G^\fg,\
C_{G,\fg}//G=\fZ_\fg^G,\ C_{G,G}//G=\fZ_G^G$. The projections to the second
(regular) factor are denoted by
$\varpi:\ \fZ_\fg^\fg\to\fg^{reg}//G=\ft/W,\
\varpi:\ \fZ_G^\fg\to G^{reg}//G=T/W,\
\varpi:\ \fZ_\fg^G\to\fg^{reg}//G=\ft/W,\
\varpi:\ \fZ_G^G\to G^{reg}//G=T/W$. In all the four cases $\varpi$ is flat.

We consider the restrictions of our centralizer varieties to the Kostant
slices and Steinberg cross-sections: $C_{\fg,\fg}^\Sigma=C_{\fg,\fg}\cap(\fg\times\Sigma_\fg),\
C_{\fg,G}^\Sigma=C_{\fg,G}\cap(\fg\times\Sigma_G),\
C_{G,\fg}^\Sigma=C_{G,\fg}\cap(G\times\Sigma_\fg),\
C_{G,G}^\Sigma=C_{G,G}\cap(G\times\Sigma_G)$.

Then the locally closed embedding $C_{\fg,\fg}^\Sigma\hookrightarrow C_{\fg,\fg}
\twoheadrightarrow\fZ_\fg^\fg$ induces an isomorphism
$C_{\fg,\fg}^\Sigma\simeq\fZ_\fg^\fg$. Similarly, we have isomorphisms
$C_{G,\fg}^\Sigma\simeq\fZ_\fg^G$ and (for simply connected $G$)
$C_{\fg,G}^\Sigma\simeq\fZ_G^\fg,\ C_{G,G}^\Sigma\simeq\fZ_G^G$. The universal centralizer $\fZ_G^G$
for simply connected $G$ was introduced by G.~Lusztig in~section~8 of~\cite{lus}.

Thus both $\fZ_\fg^\fg\to\ft/W$ and $\fZ_G^\fg\to T/W$
(for simply connected $G$) are the sheaves of abelian Lie
algebras, while both $\fZ_\fg^G\to\ft/W$ and $\fZ_G^G\to T/W$
(for simply connected $G$) are the sheaves of abelian Lie groups.

\subsection{Isogenies}
The center $Z(G)$ acts naturally on $\fZ_G^\fg$ (resp. $\fZ_\fg^G$) by
$z(x,g)=(x,zg)$ (resp. $z(g,x)=(zg,x)$). The center $Z(G)$ acts on
$\fZ_G^G$ on both sides: $z_1(g_1,g_2)z_2=(z_1g_1,z_2g_2)$.
Let $\tilde G$ denote the universal cover of $G$. Then the fundamental
group $\pi_1(G)$ is embedded into $Z({\tilde G})$, and we have
$\fZ_G^\fg=\pi_1(G)\backslash\fZ_{\tilde G}^\fg,\
\fZ_\fg^G=\pi_1(G)\backslash\fZ_\fg^{\tilde G},\
\fZ_G^G=\pi_1(G)\backslash\fZ_{\tilde G}^{\tilde G}/\pi_1(G)$.

\subsection{Symplectic structures}
\label{sym}
We fix an invariant identification $\fg\simeq\fg^*$, hence $\ft\simeq\ft^*$.
Then $\fg\times\fg$ gets identified with $\fg\times\fg^*=\sT^*\fg$
(the cotangent bundle), and $G\times\fg$ gets identified with
$G\times\fg^*=\sT^*G$. After this $\fZ_\fg^\fg$ (resp. $\fZ_\fg^G$) can be
viewed as a hamiltonian reduction of $\sT^*\fg$ (resp. $\sT^*G$); thus it
inherits a canonical symplectic structure.

Identifying $\fg\times G$ with $\fg^*\times G=\sT^*G$ we can view
$\fZ_G^\fg$ as a hamiltonian reduction of $\sT^*G$ as well;
thus it inherits a canonical Poisson structure.
Note that $\fZ_G^\fg$ is smooth and symplectic iff $G$ is simply connected.
We have symplectic isomorphisms $\fZ_\fg^\fg\simeq\sT^*(\ft/W)$, and
(in case $G$ is simply connected) $\fZ_G^\fg\simeq\sT^*(T/W)$.

Note that $\fZ_G^\fg$ and $\fZ_\fg^G$ share a common open piece ${\mathsf Z}(G)$
formed by the classes of pairs $(g,x)$ where both $g$ and $x$ are regular.
The canonical symplectic structures agree on
$\fZ_G^\fg\supset{\mathsf Z}(G)\subset\fZ_\fg^G$.
Note also that for adjoint $G$ the space ${\mathsf Z}(G)$
contains (a complexification of)
the Kostant's phase space $\fZ(G)$ of the classical Toda lattice ~\cite{k},
and the embedding $\fZ(G)\hookrightarrow\fZ_\fg^G$ is given by the Theorem ~2.6
of {\em loc. cit.}

A.~Alexeev, A.~Malkin and E.~Meinrenken introduced in ~\cite{amm} ~Example ~6.1
the q-Hamiltonian $G$-space {\em internal fusion double} ${\mathbf D}(G)$.
We have a natural map $\fZ^G_G\to{\mathbf D}(G)//\Delta_G$ (q-Hamiltonian reduction with
respect to the diagonal action of $G$). This map is a birational isomorphism, but not an
isomorphism. For example, it contracts the neutral connected component of the centralizer of
a regular unipotent element.
The q-Hamiltonian reduction ${\mathbf D}(G)//\Delta_G$ inherits a canonical
Poisson structure. Its pullback to $\fZ^G_G$ gives a rational Poisson structure on $\fZ^G_G$.
But if $G$ is simply connected, $\fZ^G_G$ is smooth, and this rational Poisson structure arises
from a (regular) symplectic structure on $\fZ^G_G$.

\subsection{Affine blow-ups}
\label{Aff}
The set of roots of $G$ (resp. $\LG$) is denoted by $R$ (resp. $\LR$).
We will view $\alpha\in R$ (resp. $\alphach\in\LR$) as a homomorphism
$\ft\to\BC$ (resp. $\Lt\to\BC$) or as a homomorphism $T\to\BC^*$
(resp. $\LT\to\BC^*$) depending on a context. Also, for a root $\alpha\in R$
we denote by $^1\alpha$ (resp. $^2\alpha$) the linear function on
$\ft\times\ft$ obtained as a composition of $\alpha$ with the
projection to the first (resp. second) factor. Finally, $\ft_\alpha\subset\ft$, $T_\alpha\subset T$, and
$\LT_{\check\alpha}\subset\LT$ denote the kernels of the corresponding root homomorphisms.

Let $D_\ft\subset \ft$, $D_T\subset T$, $D_{\LT} \subset \LT$ be the union of kernels of root homomorphisms
over all the positive roots (the discriminant). (Its ideal is the intersection of the ideals of the kernels of the corresponding root
homomorphisms.) Also, let $D_{\ft,\ft}\subset\ft\times\ft$
(resp.\ $D_{\ft,T}\subset\ft\times T$, $D_{T,\ft}\subset T\times\ft$, $D_{T,T}\subset T\times T$,
$D_{\LT,T}\subset\LT\times T$) denote the union of $\ft_\alpha\times\ft_\alpha$ (resp.\ $\ft_\alpha\times T_\alpha$,
$T_\alpha\times\ft_\alpha$, $T_\alpha\times T_\alpha$, $\LT_{\check\alpha}\times T_\alpha$) over all the positive roots.

Consider the morphisms $\fBt_\fg^\fg\to \ft\times\ft$, $\fBt^\fg_G\to \ft\times T$, $\fBt_\fg^G \to T\times \ft$,
 $\fBt_G^G\to T\times T$ and $\fBt_G^\LG\to \LT\times T$ defined as follows.
We let $\fBt_\fg^\fg$ be the open subscheme in the blow up  $Bl_{D_{\ft,\ft}}(\ft\times \ft)$ of $\ft\times \ft$ at
$D_{\ft,\ft}$ whose complement is the strict transform of the divisor $\ft\times D_\ft$. Likewise,
 $\fBt^\fg_G\subset Bl_{D_{\ft,T}}(\ft\times T)$, 
 $\fBt_\fg^G\subset  Bl_{D_{T,\ft}}(T\times \ft)$,
 $\fBt_G^G\subset  Bl_{D_{T,T}}(T\times T)$ and
   $\fBt_G^\LG\subset Bl_{D_{\LT,T}}(\LT\times T)$ are open complements
   of the strict transforms of $\ft\times D_T$, $T\times D_{\ft}$, $T\times D_T$ and $\LT \times D_T$ respectively.

We also set
$\fB_\fg^\fg=\fBt_\fg^\fg/W$, $\fB_G^G=\fBt_G^G/W$, $\fB_G^\fg=\fBt_G^\fg/W$,
$\fB^G_\fg=\fBt^G_\fg/W$, $\fB_G^\LG=\fBt_G^\LG/W$
We denote by $\varpit$, $\varpi$ the projection of $\fBt$ (resp.\ $\fB$)
to the second factor; thus we have $\varpit\colon \fBt_\fg^\fg\to\ft,\ \fBt_G^\fg\to T,\
\fBt_\fg^G\to\ft,\ \fBt_G^G\to T,\ \fB_G^\LG\to T$ and similarly for $\varpi$.

\subsection{Poisson structures}
\label{symp}
We have the canonical trivializations of the tangent bundles
$\sT(\ft\times\ft)=(\ft\times\ft)\times(\ft\times\ft),\
\sT(\ft\times T)=(\ft\times T)\times(\ft\times\ft),\
\sT(T\times\ft)=(T\times\ft)\times(\ft\times\ft),\
\sT(T\times T)=(T\times T)\times(\ft\times\ft),\
\sT(T\times\LT)=(T\times\LT)\times(\ft\times\Lt)$.
Making use of the identification $\Lt=\ft^*\simeq\ft$
we obtain the $W$-invariant
symplectic structures on the above varieties. Thus the above affine blow-ups
carry the rational Poisson structures (regular off the discriminants
$\bD\subset\fB$).

\begin{prop}
\label{simp}
The Poisson structure on $\fB_\fg^\fg-\bD$ (resp. $\fB_G^\fg-\bD,\
\fB_\fg^G-\bD,\ \fB_G^G-\bD,\ \fB_G^\LG-\bD$) extends to the
global Poisson structure;
it is a symplectic structure if the corresponding variety is smooth.
\end{prop}

\begin{prop}
\label{razdutija}
We are in the setup of ~\ref{Aff}.

a) The morphism $\varpit$ is smooth. The morphism $\varpi$ is flat if $G$ is
simply connected.

b) There are natural identifications $\fB_\fg^\fg\simeq\fZ_\fg^\fg,\
\fB_G^\fg\simeq\fZ_G^\fg,\ \fB_\fg^G\simeq\fZ_\fg^G,\ \fB_G^G\simeq\fZ_G^G$
commuting with $\varpi$.

c) If $G$ is simply laced and adjoint, we have an identification
$\fB^G_\LG\simeq Z(\LG)\backslash\fZ^{\check G}_\LG$ commuting with $\varpi$.

d) If $G$ is simply laced and simply connected, we have an identification
$\fB^G_\LG\simeq\fZ^G_G/Z(G)$ commuting with $\varpi$.

e) The above identifications respect the Poisson structures.

\end{prop}

\begin{Rem} Proposition \ref{razdutija}a) implies that varieties $\fB_\fg^\fg$, $\fBt_\fg^\fg$ etc. are smooth.
If $\g={\mathfrak{gl}}(n)$ then a deep theorem of M.~Haiman \cite{Ha} shows that the quotient of the blow up of the discriminant ideal $Bl_{D_{\ft,\ft}}(\ft\times \ft)/W$ is smooth, it is in fact identified with the Hilbert scheme of points on the plane. The analogue of Haiman's theorem is known to fail for groups of types other than $A_n$; however, Proposition \ref{razdutija}a) exhibits an open subset in the blow up that is shown to be smooth by a much simpler argument below. 
\end{Rem}

\subsection{Flat group sheaves}
We consider the functor $\fF_\fg^\fg$
on the category $\on{Flat}_{\ft/W}$ of schemes flat over
$\ft/W$ to the category of sets, sending a test scheme $S$ to
the set of $W$-invariant morphisms
$\left(\on{Mor}(S\times_{\ft/W}\ft,\ft)\right)^W$.
Similarly, we consider the functor $\fF_G^\fg$
on the category $\on{Flat}_{T/W}$
sending a test scheme $S$ to the set of $W$-invariant morphisms
$\left(\on{Mor}(S\times_{T/W}T,\ft)\right)^W$.
Also, we consider the functor $\fF_\fg^G$ on the category $\on{Flat}_{\ft/W}$
sending a test scheme $S$ to the set of $W$-invariant morphisms
$\left(\on{Mor}(S\times_{\ft/W}\ft,T)\right)_0^W\subset
\left(\on{Mor}(S\times_{\ft/W}\ft,T)\right)^W$
subject to the condition (cf. ~\cite{dg} ~4.2)
\begin{equation}
\label{Denis1}
\alpha\left(f(\alpha^{-1}(0))\right)=1\ \forall\ \alpha\in R.
\end{equation}
(note that the $W$-invariance condition automatically implies
$\alpha\left(f(\alpha^{-1}(0))\right)=\pm1\ \forall ~\alpha\in R.$)

Furthermore, we consider the functor $\fF_G^G$ on the category $\on{Flat}_{T/W}$
sending a test scheme $S$ to the set of $W$-invariant morphisms
$\left(\on{Mor}(S\times_{T/W}T,T)\right)_0^W\subset
\left(\on{Mor}(S\times_{T/W}T,T)\right)^W$
subject to the condition
\begin{equation}
\label{Denis2}
\alpha\left(f(\alpha^{-1}(1))\right)=1\ \forall\ \alpha\in R.
\end{equation}
(note that the $W$-invariance condition automatically implies
$\alpha\left(f(\alpha^{-1}(1))\right)=\pm1\ \forall ~\alpha\in R.$)

Finally, we consider the functor $\fF_G^\LG$ on the category $\on{Flat}_{T/W}$
sending a test scheme $S$ to the set of $W$-invariant morphisms
$\left(\on{Mor}(S\times_{T/W}T,\LT)\right)_0^W\subset
\left(\on{Mor}(S\times_{T/W}T,\LT)\right)^W$
subject to the condition
\begin{equation}
\label{Denis3}
\alphach\left(f(\alpha^{-1}(1))\right)=1\ \forall\ \alpha\in R.
\end{equation}
(note that the $W$-invariance condition automatically implies
$\alphach\left(f(\alpha^{-1}(1))\right)=\pm1\ \forall ~\alpha\in R.$)

The following Proposition is a generalization of ~\cite{dg} ~11.6.

\begin{prop}
\label{sechenija} Assume that $G$ is simply connected. 
The functor $\fF_\fg^\fg$ (resp. $\fF_G^\fg,\ \fF_\fg^G,\ \fF_G^G,\ \fF_G^\LG$)
is representable by the scheme $\fB_\fg^\fg$ (resp.
$\fB_G^\fg,\ \fB_\fg^G,\ \fB_G^G,\ \fB_G^\LG$).

\end{prop}

\subsection{Equivariant Borel-Moore Homology}
For the definition of convolution in equivariant Borel-Moore Homology
we refer the reader to ~\cite{cg} ~2.7, ~8.3 or ~\cite{l} Chapter 2.

We have $H^\GO_\bullet(pt)=H_\GO^\bullet(pt)=\BC[\ft/W]$, and
$H^{\GO\ltimes\Gm}_\bullet(pt)=H_{\GO\ltimes\Gm}^\bullet(pt)=
\BC[\ft/W][\hbar]$ where $\hbar$ is the generator of $H^2_{\Gm}(pt)$.
We will consider the $\BC[\ft/W]$-algebra
(resp. $\BC[\ft/W][\hbar]$-algebra) (with respect to convolution)
$H^\GO_\bullet(\Gr_G)$ (resp. $H^{\GO\ltimes\Gm}_\bullet(\Gr_G)$).
Note that setting $\hbar=0$ in
$H^{\GO\ltimes\Gm}_\bullet(\Gr_G)$ we obtain $H^\GO_\bullet(\Gr_G)$;
indeed for any group $H$, a space $X$ with an $H\times \Gm$ action, and
an $H\times \Gm$-equivariant complex $\F$ on $X$ we have a long exact sequence
$\dots \to
H^{i-2}_{H\times \Gm}(X,\F)\overset{\hbar}{\To}  H^{i}_{H\times \Gm}(X,\F)
\to H^i_H(X,\F)\to H^{i-1}_{H\times \Gm}(X,\F)\to \dots$
coming from the principal $\Gm$-bundle $E(H\times\Gm)\times_H
X \to E(H\times\Gm)\times_{H\times\Gm} X$;
if the space of $H\times\Gm$-equivariant cohomology is $\hbar$-torsion
free, then we get $H^\bu_{H}(X,\F)=H^\bu(X,\F)|_{\hbar =0}$.

\begin{thm}
\label{HBM}
a) The algebra $H^\GO_\bullet(\Gr_G)$ is commutative;

b) Its spectrum together with the projection onto $\ft/W=\Lt/W$ is
naturally isomorphic to $\fZ^{\check G}_\Lg\stackrel{\varpi}{\to}\Lt/W$;

c) The Poisson structure on $H^\GO_\bullet(\Gr_G)$ arising from the
$\hbar$-deformation $H^{\GO\ltimes\Gm}_\bullet(\Gr_G)$,
corresponds under the above identification to the
Poisson structure of ~\ref{sym} on $\fZ^{\check G}_\Lg$.

\end{thm}

\begin{rem}\label{relation_with_ginzburg}
{\em The equivariant cohomology ring
$H^\bullet_\GO(\Gr_G,\BC)=H^\bullet_\GO(\Gr_G)$ was computed by
V.~Ginzburg ~\cite{g}.
More precisely, 
the projection to the
second (regular) factor $\fZ^{\check\fg}_\Lg\to\Lg{}^{reg}//\LG=\Lt/W$
makes $\fZ^{\check\fg}_\Lg$ a sheaf of abelian Lie algebras. V.~Ginzburg
identifies $H^\bullet_\GO(\Gr_G)$ with the global sections of
the relative universal enveloping algebra
$U_{\Lt/W}\left(\fZ^{\check\fg}_\Lg\right)$. One can easily
check that this result is compatible with our Theorem \ref{HBM}(b) as follows.
For a group scheme $A$ over a base $S$ one has a natural pairing
$ U({\mathfrak a})\times \CO(A)\to \CO(S)$ where $U({\mathfrak a})$
is the enveloping (over $\O(S)$) of  the Lie algebra of $A$;
the pairing sends $(\xi,f)$ to $\xi(f)$ restricted to the identity of $A$.
On the other hand, for a compact (or ind-compact)
 $H$-space $X$ we have a pairing
$H^\bullet_H(X)\times H_\bullet^H(X)\to H^\bullet_H(pt)$ induced
by the action of cohomology on homology, and the push-forward 
map in Borel-Moore homology $ H_\bullet^H(X)\to H^\bullet_H(pt)$.
The isomorphisms of \cite{g} and of Theorem ~\ref{HBM} take the first
pairing into the second one.}
\end{rem}

\subsection{Equivariant $K$-theory}
For the definition of convolution in equivariant $K$-theory we refer the
reader to Chapter 5 of ~\cite{cg}.

We have $K^\GO(pt)=\BC[T/W]$, and
$K^{\GO\ltimes\Gm}(pt)=\BC[T/W][q^{\pm1}]$.
We will consider the $\BC[T/W]$-algebra
(resp. $\BC[T/W][q^{\pm1}]$-algebra) (with respect to convolution)
$K^\GO(\Gr_G)$ (resp. $K^{\GO\ltimes\Gm}(\Gr_G)$).
Note that setting $q=1$ in
$K^{\GO\ltimes\Gm}(\Gr_G)$ we obtain $K^\GO(\Gr_G)$.

\begin{thm}
\label{K}
a) The algebra $K^\GO(\Gr_G)$ is commutative;

b) Its spectrum together with the projection onto $T/W$ is
naturally isomorphic to $\fB^{\check G}_G\stackrel{\varpi}{\to}T/W$;

c) The Poisson structure on $K^\GO(\Gr_G)$ arising from the
$q$-deformation $K^{\GO\ltimes\Gm}(\Gr_G)$,
corresponds under the above identification to the
Poisson structure of ~\ref{simp} on $\fB^{\check G}_G$
in case the latter variety is smooth, i.e. $G$ is simply connected.
\end{thm}

\section{Calculations in rank 1}
\label{rank1}

In this section $G=SL_2$, and $\LG=PGL_2$.
The Weyl group $W=\BZ/2\BZ$, the Cartan torus $T=\Gm=\BC^*$ with
a coordinate $z$, and the only simple root $\alpha(z)=z^2$.
The dual torus $\LT=\Gm=\BC^*$ with a coordinate $t$, and $\alphach(t)=t$.
The Cartan Lie algebra $\ft=\BC$ with a coordinate $x=\alpha(x)$.
We fix a $\sqrt{-1}$.

\subsection{$\fZ_G^G$ and $\fB_G^G$}
\label{gruppa-gruppa}
We choose the standard $\fsl_2$-triple
$e=\begin{pmatrix}0&1\\0&0\end{pmatrix},\
h=\begin{pmatrix}1&0\\0&-1\end{pmatrix},\
f=\begin{pmatrix}0&0\\1&0\end{pmatrix}$. The standard choice of a Steinberg cross section is
$\Sigma_G=\left\{\begin{pmatrix}a&1\\-1&0\end{pmatrix},\ a\in\BC\right\}$.
However, for historical reasons we conjugate it to the following one we will work with from now on:
$\Sigma_G=\left\{\begin{pmatrix}a-1&a-2\\1&1\end{pmatrix},\ a\in\BC\right\}$.
One checks that a matrix
$\begin{pmatrix}x_{11}&x_{12}\\x_{21}&x_{22}\end{pmatrix}$ commutes with
$\begin{pmatrix}a-1&a-2\\1&1\end{pmatrix}$ iff
$\begin{pmatrix}x_{11}&x_{12}\\x_{21}&x_{22}\end{pmatrix}=
\sqrt{-1}\begin{pmatrix}(1-a)c+b&(2-a)c\\-c&b-c\end{pmatrix}$ for $b,c\in\BC$.
Then the condition $\det=1$ reads as

\begin{equation}
\label{star decomposition}
1=abc-b^2-c^2.
\end{equation}

Thus, $\fZ_G^G$ is identified with a hypersurface $\CS$ in $\BA^3$
given by the equation ~(\ref{star decomposition}).
The left (resp. right) multiplication by
$-1\in Z(SL_2)$ is an involution $\imath$ (resp. $\jmath$) on $\CS$ given by
$\imath(a,b,c)=(a,-b,-c)$ (resp. $\jmath(a,b,c)=(-a,b,-c)$).
Hence, $\fZ_{\check G}^\LG=\imath\backslash\CS/\jmath$.

Generically, we can diagonalize two commuting matrices simultaneously, that
is, there is $g\in SL_2$ such that
$g\sqrt{-1}\begin{pmatrix}(1-a)c+b&(2-a)c\\-c&b-c\end{pmatrix}g^{-1}=
\begin{pmatrix}y&0\\0&y^{-1}\end{pmatrix}$ and
$g\begin{pmatrix}a-1&a-2\\1&1\end{pmatrix}g^{-1}=
\begin{pmatrix}z&0\\0&z^{-1}\end{pmatrix}$ for some $y,z\in\Gm=\BC^*=T$
defined up to simultaneous inversion.
Then we have
\begin{equation}
\label{Moka}
a=z+z^{-1},\
b=\frac{-\sqrt{-1}}{2}
\left(y+y^{-1}+\frac{(y-y^{-1})(z+z^{-1})}{z-z^{-1}}\right),\
c=-\sqrt{-1}\frac{y-y^{-1}}{z-z^{-1}}.
\end{equation}
We conclude that $\BC[\CS]=\BC[y^{\pm1},z^{\pm1},\frac{y-y^{-1}}{z-z^{-1}}]^W$
where the nontrivial element $w\in W$ acts by $w(y,z)=(y^{-1},z^{-1})$.
We can rewrite $\BC[y^{\pm1},z^{\pm1},\frac{y-y^{-1}}{z-z^{-1}}]^W$ as
$\BC[y^{\pm1},z^{\pm1},\frac{y^2-1}{z^2-1}]^W$ to manifest its coincidence
with $\BC[\fB_G^G]$. All in all, we have $\fB_G^G\simeq\CS\simeq\fZ_G^G$.
Since we can identify $\LT$ with $T/Z(G)$, the
identifications $\fB^G_\LG\simeq\CS/\jmath,\
\fB_G^\LG\simeq\imath\backslash\CS,\
\fB_{\check G}^\LG\simeq\imath\backslash\CS/\jmath\simeq\fZ_{\check G}^\LG$
follow immediately.

\subsection{$\fZ_\fg^G$ and $\fB_\fg^G$}
\label{gruppa-algebra}
The Kostant slice
$\Sigma_\fg=\left\{\begin{pmatrix}0&\delta\\1&0\end{pmatrix},\
\delta\in\BC\right\}$.
One checks that a matrix
$\begin{pmatrix}x_{11}&x_{12}\\x_{21}&x_{22}\end{pmatrix}$ commutes with
$\begin{pmatrix}0&\delta\\1&0\end{pmatrix}$ iff
$\begin{pmatrix}x_{11}&x_{12}\\x_{21}&x_{22}\end{pmatrix}=
\begin{pmatrix}\xi&\delta\eta\\ \eta&\xi\end{pmatrix}$ for $\xi,\eta\in\BC$.
Then the condition $\det=1$ reads as

\begin{equation}
\label{det=1}
1=\xi^2-\delta\eta^2.
\end{equation}

Thus, $\fZ_\fg^G$ is identified with a hypersurface $\CS'$ in $\BA^3$
given by the equation ~(\ref{det=1}).
The action of
$-1\in Z(SL_2)$ is an involution $\imath$ on $\CS'$ given by
$\imath(\delta,\xi,\eta)=(\delta,-\xi,-\eta)$.
Hence, $\fZ^{\check G}_\Lg=\imath\backslash\CS'$.

Generically, we can diagonalize two commuting matrices simultaneously, that
is, there is $g\in SL_2$ such that
$g\begin{pmatrix}\xi&\delta\eta\\ \eta&\xi\end{pmatrix}g^{-1}=
\begin{pmatrix}y&0\\0&y^{-1}\end{pmatrix}$ and
$g\begin{pmatrix}0&\delta\\1&0\end{pmatrix}g^{-1}=
\begin{pmatrix}x&0\\0&-x\end{pmatrix}$ for some $y\in\Gm=\BC^*=T,\
x\in\BC=\ft$, defined up to $(y,x)\mapsto(y^{-1},-x)$.
Then we have $$\delta=x^2,\
\xi=\frac{y+y^{-1}}{2},\
\eta=\frac{y-y^{-1}}{2x}.$$
We conclude that $\BC[\CS']=\BC[y^{\pm1},x,\frac{y-y^{-1}}{x}]^W$
where the nontrivial element $w\in W$ acts by $w(y,x)=(y^{-1},-x)$.
We can rewrite $\BC[y^{\pm1},x,\frac{y-y^{-1}}{x}]^W$ as
$\BC[y^{\pm1},x,\frac{y^2-1}{x}]^W$ to manifest its coincidence
with $\BC[\fB_\fg^G]$. All in all, we have $\fB_\fg^G\simeq\CS'\simeq\fZ_\fg^G$.
Since we can identify $\LT$ with $T/Z(G)$, the identification
$\fB^{\check G}_\Lg\simeq\imath\backslash\CS'\simeq\fZ^{\check G}_\Lg$
follows immediately.

\subsection{$\fZ_G^\fg$ and $\fB_G^\fg$}
\label{algebra-gruppa}
Recall the Steinberg cross-section
$\Sigma_G=\left\{\begin{pmatrix}a-1&a-2\\1&1\end{pmatrix},\ a\in\BC\right\}$.
One checks that a traceless matrix
$\begin{pmatrix}x_{11}&x_{12}\\x_{21}&-x_{11}\end{pmatrix}$ commutes with
$\begin{pmatrix}a-1&a-2\\1&1\end{pmatrix}$ iff
$\begin{pmatrix}x_{11}&x_{12}\\x_{21}&-x_{11}\end{pmatrix}=
\zeta\begin{pmatrix}2-a&4-2a\\-2&a-2\end{pmatrix}$ for $\zeta\in\BC$.

Thus, $\fZ_G^\fg$ is identified with $\BA^2$ with coordinates $a,\zeta$.
The action of
$-1\in Z(SL_2)$ is an involution $\jmath$ on $\BA^2$ given by
$\jmath(a,\zeta)=(-a,-\zeta)$.
Hence, $\fZ_{\check G}^\Lg=\BA^2/\jmath$.

Generically, we can diagonalize two commuting matrices simultaneously, that
is, there is $g\in SL_2$ such that
$g\zeta\begin{pmatrix}2-a&4-2a\\-2&a-2\end{pmatrix}g^{-1}=
\begin{pmatrix}x&0\\0&-x\end{pmatrix}$ and
$g\begin{pmatrix}a-1&a-2\\1&1\end{pmatrix}g^{-1}=
\begin{pmatrix}z&0\\0&z^{-1}\end{pmatrix}$ for some $x\in\BC=\ft,\
z\in\Gm=\BC^*=T$ defined up to $(x,z)\mapsto(-x,z^{-1})$.
Then we have $$a=z+z^{-1},\
\zeta=\frac{x}{z-z^{-1}}.$$
We conclude that $\BC[\BA^2]=\BC[x,z^{\pm1},\frac{x}{z-z^{-1}}]^W$
where the nontrivial element $w\in W$ acts by $w(x,z)=(-x,z^{-1})$.
We can rewrite $\BC[x,z^{\pm1},\frac{x}{z-z^{-1}}]^W$ as
$\BC[x,z^{\pm1},\frac{x}{z^2-1}]^W$ to manifest its coincidence
with $\BC[\fB_G^\fg]$.
All in all, we have $\fB_G^\fg\simeq\BA^2\simeq\fZ_G^\fg$.
Since we can identify $\LT$ with $T/Z(G)$, the
identification $\fB_{\check G}^\Lg\simeq\BA^2/\jmath\simeq\fZ_{\check G}^\Lg$
follows immediately.

\subsection{$\fZ_\fg^\fg$ and $\fB_\fg^\fg$}
\label{algebra-algebra}
Recall the Kostant slice
$\Sigma_\fg=\left\{\begin{pmatrix}0&\delta\\1&0\end{pmatrix},\
\delta\in\BC\right\}$.
One checks that a traceless matrix
$\begin{pmatrix}x_{11}&x_{12}\\x_{21}&-x_{11}\end{pmatrix}$ commutes with
$\begin{pmatrix}0&\delta\\1&0\end{pmatrix}$ iff
$\begin{pmatrix}x_{11}&x_{12}\\x_{21}&-x_{11}\end{pmatrix}=
\begin{pmatrix}0&\delta\theta\\ \theta&0\end{pmatrix}$ for $\theta\in\BC$.
Thus, $\fZ_\fg^\fg$ is identified with $\BA^2$ with coordinates $\delta,\theta$.

Generically, we can diagonalize two commuting matrices simultaneously, that
is, there is $g\in SL_2$ such that
$g\begin{pmatrix}0&\delta\theta\\ \theta&0\end{pmatrix}g^{-1}=
\begin{pmatrix}u&0\\0&-u\end{pmatrix}$ and
$g\begin{pmatrix}0&\delta\\1&0\end{pmatrix}g^{-1}=
\begin{pmatrix}x&0\\0&-x\end{pmatrix}$ for some
$u,x\in\BC=\ft$, defined up to $(u,x)\mapsto(-u,-x)$.
Then we have $$\delta=x^2,\
\theta=\frac{u}{x}.$$
We conclude that $\BC[\BA^2]=\BC[u,x,\frac{u}{x}]^W$
where the nontrivial element $w\in W$ acts by $w(u,x)=(-u,-x)$.
Hence we get an identification $\fB_\fg^\fg\simeq\BA^2\simeq\fZ_\fg^\fg$.

\subsection{$\fB_\fg^G$ and $\fF_\fg^G$}
Recall the setup of Proposition ~\ref{sechenija}. We will prove that the
functor $\fF_\fg^G$ is representable by the scheme $\fB_\fg^G$; the other parts
of the Proposition are proved absolutely similarly, as well as the Proposition
for $G$ replaced by $\LG$. For a scheme $S$ flat over
$\ft/W$ we will denote by $S_\ft$ the cartesian product $S\times_{\ft/W}\ft$.
Our usual coordinate $x$ on $\ft$
gives rise to the same named function on $S_\ft$.
The nontrivial element $w\in W$ acts by the involution of $S_\ft$.
Finally, we denote by $(\fB_\fg^G)_{S_\ft}$
the affine blow-up of $S_\ft\times T$, that is $S_\ft\times_\ft\fB_\fg^G$.
Clearly, $w$ acts as an involution of $(\fB_\fg^G)_{S_\ft}$.

Note that the condition ~(\ref{Denis1}) is void in the case under consideration.
Given a $w$-equivariant morphism $f:\ S_\ft\to T=\Gm$ we see that $f^2-1$
is divisible by $x$, hence $f$ lifts uniquely to a section $\hat f$ of
$(\fB_\fg^G)_{S_\ft}$ over $S_\ft$. Evidently, $\hat f$ is $w$-invariant.
If we consider $\hat f$ as a closed subscheme of $(\fB_\fg^G)_{S_\ft}$,
then ${\hat f}/W$ is a closed subscheme of
$(\fB_\fg^G)_{S_\ft}/W=S\times_{\ft/W}\fB_\fg^G$ which is the graph of a
morphism ${\tilde f}:\ S\to\fB_\fg^G$.

Conversely, given a morphism ${\tilde f}:\ S\to\fB_\fg^G$ we consider its
graph $\Gamma_{\tilde f}$ as a closed subscheme of $S\times_{\ft/W}\fB_\fg^G$,
and then the cartesian product
$\Gamma_{\tilde f}\times_{S\times_{\ft/W}\fB_\fg^G}(\fB_\fg^G)_{S_\ft}$
is a section $\hat f$ of $(\fB_\fg^G)_{S_\ft}$ over $S_\ft$.
Evidently, $\hat f$ gives rise to a $w$-equivariant function $f:\ S_\ft\to T$.

\medskip




\subsection{A basis in equivariant $K$-theory}
\label{K-teorija}
We recall a few standard facts about the affine Grassmannians $\Gr_G$ and
$\Gr_\LG$. The $\GO$-orbits (equivalently, $\LGO$-orbits) on $\Gr_\LG$
are numbered by nonnegative integers and denoted by $\Gr_{\LG,n},\ n\in\BN$.
The orbits $\Gr_{\LG,2n},\ n\in\BN$, form a connected component of $\Gr_\LG$
equal to $\Gr_G$. The open embedding of an orbit into its closure will be
denoted by $j_n:\ \Gr_{\LG,n}\hookrightarrow\ol{\Gr}_{\LG,n}$ or simply by $j$
if no confusion is likely. We have $\dim\Gr_{\LG,n}=n$; in particular,
$\Gr_{\LG,0}$ is a point.

We have $K^\GO(\Gr_{\LG,0})=\on{Rep}(G)$ with a basis $\bv(n),\ n\in\BN$,
formed by the classes of irreducible $G$-modules $\CV(n)$.
Also, $K^\LGO(\Gr_{\LG,0})=\on{Rep}(\LG)\subset\on{Rep}(G)$ has a basis
$\bv({2n}),\ n\in\BN$.

For $m>0$ the $\GO$-equivariant line bundles in $\Gr_{\LG,m}$ are numbered
by integers and denoted by $\CL(n)_m$. Among them, the $\LGO$-equivariant
line bundles are exactly $\CL(2n)_m,\ n\in\BZ$.
We define $\CV(n)_m$ as $j_*\CL(n)_m[\frac{m}{2}]$, that is, the
(nonderived) direct image to the orbit closure placed in the homological degree
$-\frac{m}{2}$. Note that since the complement $\ol{\Gr}_{\LG,m}-\Gr_{\LG,m}$
has codimension 2, the above direct image is a coherent sheaf. The degree shift
will become clear later.
The class $[\CL(n)_m]$ in $K^\GO(\Gr_\LG)$ will be denoted by $\bv(n)_m$.
Thus, it is natural to denote $\bv(n)$ above by $\bv(n)_0$; we will keep
both names.

The collection $\{\bv(n)_m:\ n\in\BN$ if $m=0;\ n\in\BZ$ if $m\in\BN-0\}$
forms a basis in $K^\GO(\Gr_\LG)$. Among this collection, all the
$\bv(n)_m$ with $n$ even (resp. $m$ even) form a basis in
$K^\LGO(\Gr_\LG)$ (resp. $K^\GO(\Gr_G)$).

\subsection{Convolution: commutativity}
\label{kommut}
In this subsection $G$ is an arbitrary semisimple group. We prove ~\ref{K} ~(a).
We refer the reader to ~\cite{ga} for the basics of Beilinson-Drinfeld
Grassmannian. Recall that $\Gr^{BD}_G\stackrel{\pi}{\to}\BA^1$
is a flat ind-scheme such that
$\pi^{-1}(\BA^1-0)=(\BA^1-0)\times\Gr_G\times\Gr_G$, while
$\pi^{-1}(0)=\Gr_G$. We also have the deformed convolution diagram
$\Gr^{BD,conv}_G\stackrel{\Pi}{\to}\Gr^{BD}_G$ such that $\Pi$ is an
isomorphism over $\BA^1-0$, while over $0\in\BA^1$ our $\Pi$
is the usual convolution diagram
$\GF\times_\GO\Gr_G\stackrel{\Pi_0}{\to}\Gr_G$.

Given two $\GO$-equivariant complexes of coherent sheaves $\CA,\CB$ on $\Gr_G$,
we can form their ``deformed convolution" complex $\CA{\widetilde\star}\CB$
on $\Gr^{BD,conv}_G$ such that over $\BA^1-0$ it is isomorphic to
$\CO_{\BA^1-0}\boxtimes\CA\boxtimes\CB$, while over $0\in\BA^1$ it is
isomorphic to the usual twisted product $\CA\ltimes\CB$ on the convolution
diagram $\GF\times_\GO\Gr_G$. In addition, if $\CA,\CB$ are coherent sheaves,
then $\CA{\widetilde\star}\CB$ is flat over $\BA^1$.
It implies that in the $K$-group the class
$[\CA\ltimes\CB]$ is the {\em specialization} (see ~\cite{cg} 5.3)
of the class $[\CO_{\BA^1-0}\boxtimes\CA\boxtimes\CB]$ in the family
$\Gr^{BD,conv}_G\stackrel{\pi\circ\Pi}{\longrightarrow}\BA^1$,
and also the class
$[\CA\star\CB]=[\Pi_{0*}(\CA\ltimes\CB)]$ is the specialization
of the class $[\CO_{\BA^1-0}\boxtimes\CA\boxtimes\CB]$ in the family
$\Gr^{BD}_G\stackrel{\pi}{\longrightarrow}\BA^1$.
Hence the desired commutativity.

\subsection{Convolution: relations}
\label{sootn}
We return to the setup of ~\ref{K-teorija}. Note that $\Gr_{\LG,1}\simeq\BP^1$,
and $\CV(n)_1$ is the line bundle $\CO(n)$ on $\BP^1$.
The twisted product $\CV(n)_1\ltimes\CV(l)_1$ is the line bundle $\CO(n,l)$
on the 2-dimensional subvariety
$\CH_2\subset{\check G}({\mathbf F})\times_\LGO\Gr_\LG$ isomorphic to the
Hirzebruch surface $\BP(\CO(2)\oplus\CO)$ over $\BP^1$. The projection
$\Pi_0:\ \CH_2\to\Gr_{\LG,2}$ is the contraction of the $-2$-section
$\BP^1\hookrightarrow\CH_2$.

Now it is easy to compute $\bv(n)_1\star\bv(n)_1=\bv(2n)_2,\
\bv(1)_1\star\bv(-1)_1=\bv(0)_2+1$. Taking into account the evident
relation $\bv(1)_0\star\bv(0)_1=\bv(1)_1+\bv(-1)_1$ we arrive at

\begin{equation}
\label{Ivan}
\bv(1)_0\star\bv(0)_1\star\bv(1)_1=
\bv(1)_1\star\bv(1)_1+\bv(0)_1\star\bv(0)_1+1.
\end{equation}

A moment of reflection shows that $K^\GO(\Gr_G)$ is generated as algebra by
$\bv(1)_0,\ \bv(0)_2=\bv(0)_1\star\bv(0)_1,\ \bv(2)_2=\bv(1)_1\star\bv(1)_1,\
\bv(1)_2=\bv(1)_1\star\bv(0)_1$ (one has to use that
$\bv(k)_{2l}\star\bv(n)_{2m}=\bv(k+n)_{2l+2m}$ plus the terms supported on
the smaller orbits).
Similarly, $K^\LGO(\Gr_\LG)$ is generated as algebra by
$\bv(2)_0=\bv(1)_0\star\bv(1)_0-1,\ \bv(0)_1,\ \bv(2)_2=\bv(1)_1\star\bv(1)_1,\
\bv(2)_1=\bv(1)_1\star\bv(1)_0-\bv(0)_1$.

Note that both algebras $K^\GO(\Gr_G)$ and $K^\LGO(\Gr_\LG)$ lie in the
vector space $K^\GO(\Gr_\LG)$, and their intersection is the common
subalgebra $K^\LGO(\Gr_G)$. The tensor product algebra
$K^\GO(\Gr_G)\otimes_{K^\LGO(\Gr_G)}K^\LGO(\Gr_\LG)$ can be identified as
a vector space with $K^\GO(\Gr_\LG)$, and then it is generated by
the three basic elements $\bv(1)_0,\bv(0)_1,\bv(1)_1$ subject to the only
relation ~(\ref{Ivan}).

The comparison of equations ~(\ref{Ivan}) and ~(\ref{star decomposition})
shows that the assignment $a\mapsto\bv(1)_0,\ b\mapsto\bv(0)_1,\
c\mapsto\bv(1)_1$ establishes an isomorphism $\BC[\CS]\simeq K^\GO(\Gr_\LG)$.
It identifies the spectrum of $K^\GO(\Gr_G)$ with
$\imath\backslash\CS\simeq\fB_G^\LG$, and the spectrum of $K^\LGO(\Gr_\LG)$
with $\CS/\jmath\simeq\fB^G_\LG$.

\subsection{Iwahori-equivariant $K$-theory}
\label{Iwahori}
Let $I\subset\GO$ be the Iwahori subgroup.
The space $K^I(\Gr_G)=K^T(\Gr_G)=K^\TO(\Gr_G)=K(\TO\backslash\GF/\GO)$
is equipped with
the two commuting actions: $K(\TO\backslash\TF/\TO)$ acts by convolutions
on the left, and $K^G(\Gr_G)=K^\GO(\Gr_G)=K(\GO\backslash\GF/\GO)$
acts by convolutions on the right. Also, $W$ acts on $K^T(\Gr_G)$ commuting
with the right action of $K^G(\Gr_G)$.
Clearly, the algebra $K(\TO\backslash\TF/\TO)$ is isomorphic to
$\BC[\LT\times T]$. The action of $W$ on $K^T(\Gr_G)$ normalizes the action of
$K(\TO\backslash\TF/\TO)$ and induces the natural (diagonal) action of $W$
on $\BC[\LT\times T]$.

Our aim in this subsection is to identify the
$K(\TO\backslash\TF/\TO)\ltimes W-K^G(\Gr_G)$-bimodule $K^T(\Gr_G)$ with the
$\BC[\LT\times T]\ltimes W-\BC[\fB_G^\LG]$-bimodule $\BC[\fBt_G^\LG]$
(and similarly
for $G$ replaced by $\LG$). As in ~\ref{sootn}, it suffices to identify the
$K(\TO\backslash\LT({\mathbf F})/\LT({\mathbf O}))\ltimes
W-K^G(\Gr_\LG)$-bimodule $K^T(\Gr_\LG)$ with the
$\BC[T\times T]\ltimes W-\BC[\fB_G^G]$-bimodule $\BC[\fBt_G^G]$.

Note that $K^G(\Gr_\LG)\subset K^T(\Gr_\LG)$, and the $K^G(\Gr_\LG)$-module
$K^T(\Gr_\LG)$ is free of rank 2 with the generators $1,z$ where
$z$ is the generator of $K^T(pt)=\BC[T]$ (so that, e.g. $\bv(1)_0=z+z^{-1}$).
Furthermore, $\BC[y^{\pm1},z^{\pm1}]=\BC[T\times T]=
K(\TO\backslash\LT({\mathbf F})/\LT({\mathbf O}))\subset K^T(\Gr_\LG)$,
and one can check that


\begin{equation}
\label{moka}
y+y^{-1}=\sqrt{-1}(2\bv(0)_1-\bv(1)_0\star\bv(1)_1),\
y-y^{-1}=\sqrt{-1}(z-z^{-1})\bv(1)_1.
\end{equation}
Comparing ~(\ref{moka}) with ~(\ref{Moka}) we get the desired identification
of the $K(\TO\backslash\LT({\mathbf F})/\LT({\mathbf O}))\ltimes
W-K^G(\Gr_\LG)$-bimodule $K^T(\Gr_\LG)$ with the
$\BC[y^{\pm1},z^{\pm1}]\ltimes W-
\BC[y^{\pm1},z^{\pm1},\frac{y-y^{-1}}{z-z^{-1}}]^W$-bimodule
$\BC[y^{\pm1},z^{\pm1},\frac{y-y^{-1}}{z-z^{-1}}]$.

\subsection{Borel-Moore Homology}
\label{BMH}
For an arbitrary semisimple $G$ one proves the commutativity of
$H^\GO_\bullet(\Gr_G)$ (Theorem ~\ref{HBM} a) exactly as in ~\ref{kommut}
using the Beilinson-Drinfeld Grassmannian and the {\em specialization}
in Borel-Moore Homology (see ~\cite{cg} 2.6.30).

For $\LG=PGL_2$, let us denote by $\delta\in H^4_\LGO(pt,\BZ)=H_4^\LGO(pt,\BZ)$
the generator of the equivariant (co)homology. Furthermore, we denote by
$\eta$ (resp. $\xi$) the generator of $H_{-2}^\LGO(\Gr_{\LG,1},\BZ)$
(resp. the generator of $H_{0}^\LGO(\Gr_{\LG,1},\BZ)$).
Then it is easy to see that $\delta,\xi,\eta$ generate
$H^\LGO_\bullet(\Gr_\LG)$
(while $\delta,\xi^2,\eta^2,\xi\eta$ generate the subalgebra
$H^\GO_\bullet(\Gr_G)$),
and we claim that

\begin{equation}
\label{Roma}
1=\xi^2-\delta\eta^2.
\end{equation}

In effect, this is an equality in $H_0^\LGO(\Gr_{\LG,2})$. Since $\Gr_{\LG,2}$
is rationally smooth, $H_0^\LGO(\Gr_{\LG,2})=H^4_\LGO(\Gr_{\LG,2})$.
Let us denote by $\bB\Gr_{\LG,2}\stackrel{p}{\to}\bB\LGO$
the associated fibre bundle over the
classifying space of $\LGO$ with the fiber $\Gr_{\LG,2}$.
Then $1\in H^4_\LGO(\Gr_{\LG,2})=H^4(\bB\Gr_{\LG,2})$
is the Poincar\'e dual class of the codimension 2 cycle
$\bB\LGO=\bB\Gr_{\LG,0}\hookrightarrow\bB\Gr_{\LG,2}$,
and $\delta\eta^2=p^*\delta$.

Recall the convolution morphism
$\Pi_0:\ \CH_2\to\Gr_{\LG,2}$ of ~\ref{sootn}. This is a morphism of
$\LGO$-varieties, and we denote by $\Pi_0:\ \bB\CH_2\to\bB\Gr_{\LG,2}$
the corresponding morphism of associated fibre bundles. Note that
(additively) $H^\bullet(\bB\CH_2)=
H^\bullet(\bB\Gr_{\LG,1})\otimes_{H^\bullet(\bB\LGO)}H^\bullet(\bB\Gr_{\LG,1})$.
Recall that $\xi$ is the generator of
$H_0^\LGO(\Gr_{\LG,1})=H^2_\LGO(\Gr_{\LG,1})=H^2(\bB\Gr_{\LG,1})$.
Finally, we have $\xi^2=\Pi_{0*}(\xi\otimes\xi)$. Now ~(\ref{Roma}) follows
easily.

Comparing the sizes of $H^\LGO_\bullet(\Gr_\LG)$ and
$\BC[\delta,\xi,\eta]/(\xi^2-\delta\eta^2-1)$ we conclude that they are
isomorphic. The comparison with the equation ~(\ref{det=1}) establishes
an isomorphism $H^\LGO_\bullet(\Gr_\LG)\simeq\BC[\CS']$, and identifies
the spectrum of $H^\LGO_\bullet(\Gr_\LG)$ with $\CS'\simeq\fZ_\fg^G$,
and the spectrum of $H^\GO_\bullet(\Gr_G)$ with
$\imath\backslash\CS'\simeq\fZ^{\check G}_\Lg$.

\section{Centralizers and blow-ups}
\label{proof28}

The aim of this section is a proof of Proposition ~\ref{razdutija}.

Till the further notice $G$ is assumed simply connected. 

\subsection{Alternative descriptions of affine blow ups}

Before proceeding to the proof, we present several alternative descriptions of varieties defined
in section \ref{Aff}.

1) For a root $\alpha$ let $S_\alpha = {\mathbb C}[{\mathfrak t} \times
  {\mathfrak t}][\frac{^1\alpha}{^2\alpha},\ {}^2\beta_i^{-1}]$
where $\beta_i$ runs over all roots different from $\pm \alpha$. 
Let $S_1$ be the intersection of the rings $S_\alpha$ over all roots $\alpha$, taken inside the field of rational functions on ${\mathfrak t} \times {\mathfrak t}$,
and set $X_1=Spec(S_1)$.

2) Let $I$ be the defining ideal of $D_{\ft,\ft} \subset \ft\times \ft$ and let
$I^{(n)}$ be its symbolic powers (see~\cite[Section 3.9]{eis} or~\cite[Introduction]{ddsg}). Then let
$S_2 = \mathbb C[{\mathfrak t} \times {\mathfrak t}] [f/\Delta^n], f\in I^{(n)}$, and let
$\Delta =\prod_{\alpha\in R} {}^2\alpha$ be the discriminant pulled back from the second
factor, $X_2:=Spec(S_2)$.

3) Let $X_3$ be the open subscheme in the symbolic blow up of ${\mathfrak t} \times {\mathfrak t}$ centered at
$D_{\ft,\ft}$, whose complement is the strict transform of the divisor $\ft \times D_{\ft}$.

\begin{prop}\label{blow}
We have $X_1/W\cong X_2/W\cong X_3/W \cong \fZ_\fg^\fg$.

The parallel isomorphisms hold for $\fZ_\fg^G$, $\fZ_G^\fg$, $\fZ_G^G$.
\end{prop}

\proof The isomorphism $X_2\cong X_3$ follows from the definition of symbolic blow up.
The $n$-th symbolic power of an ideal of a reduced subscheme is the intersection of the $n$-th symbolic powers of the ideals of its
components;
for a smooth component, a local function lies in the $n$-th power of its ideal iff it does so at the generic point of the component.
This  implies that $X_1=X_2$.

The rest of the proof is based on flatness of the morphisms $\fZ_\fg^\fg\to \ft/W$,
$\fZ_\fg^G\to \ft/W$, $\fZ_G^\fg\to T/W$, $\fZ_G^G\to T/W$.

We will say that $x\in \ft$ 
 is almost regular if it lies in the 
kernel of at most one root homomorphism. Let
$\ft^\bullet\subset \ft$, $T^\bullet\subset T$ be the open set of almost regular elements.

Let $X_i^\bullet$ ($i=1,2,3$), $(\fZ_\fg^\fg)^\bullet$ be the preimages of $\ft^\bullet$ in $X_i$, $\fZ_\fg^\fg$.
Existence of  canonical isomorphisms
$$X_1^\bullet/W\cong X_2^\bullet/W\cong X_3^\bullet/W \cong (\fZ_\fg^\fg)^\bullet$$
follows from the case of $G=SL(2)$ worked out in sections ~\ref{gruppa-gruppa}--\ref{algebra-algebra}.
To reduce to that case, observe that each of the spaces $\fZ_\fg^\fg$,
$\fZ_\fg^G$, $\fZ_G^\g$, $\fZ_G^G$ can be realized as the (coarse) quotient of
the set of commuting pairs in $\fg \times \fg$ etc. such that the second element in the pair is regular modulo
the action of $G$; thus it receives a map from a similar space for a Levi subgroup $L\subset G$, which is easily seen to be \' etale
on the open set of pairs where the centralizer of the semi-simple part of the second element is contained in $L$.
If $x$, $g$ are elements of $\g$, $G$ respectively such that in the Jordan decomposition $x=x_s + x_u$, $g=su$ the semi-simple
parts $x_s$, $s$ lie in $  \ft^\bullet$, $ T^\bullet$ respectively, then the centralizer $G_s$ of $x_s$, $s$ is a Levi subgroup whose derived group is
isomorphic to $SL(2)$. The computations of  sections ~\ref{gruppa-gruppa}--\ref{algebra-algebra} show that the statement holds for $G_s$; since
the map $\fZ_{\fg_s}^{\fg_s}\to \fZ_\fg^\fg$ etc. is \' etale at a pair with second element $x$ or $g$, this implies the isomorphism with $X_i^\bullet/W$. 

On the other hand, each of the quasi-coherent sheaves $\O_{X_i}$, $\O_{\fZ_\fg^\fg}$
coincides with the direct image of its restriction to the open subvariety $\ft^\bullet$:
for $X_i$ this follows from definitions, for $\fZ_\fg^\fg$ from flatness of the universal centralizer.
The claim about $\fZ_\fg^\fg$ follows, the other cases are treated similarly. \qed






\subsection{The proof of Proposition ~\ref{razdutija}}

We first prove part b) for $\fZ_\g^\g$.

Since symbolic blow up maps to the regular blow up, in view of Proposition \ref{blow}
we have a map $\fZ_\g^\g\to \fB_\g^\g$. To construct the inverse map we establish an 
explicit isomorphism $\fZ_\g^\g\cong {\mathbb A}^{2r}$ and show that coordinate functions
on $ {\mathbb A}^{2r}$ come from regular functions on $\fB_\g^\g$. 

Let $P_1,\dots, P_r$ be a system of free generators for the ring of regular functions on $\ft/W$. 

Identifying $T^*\g$ with $\g\times \g$ we get  polynomial maps $m_i=dP_i:\g\to \g$. It is easy to see that $[x,m_i(x)]=0$ for all $x\in \g$, $i=1,\dots, r$.
The following properties of these maps can be found in (the proof of)~Lemma~6.7.30
and~Theorem~6.7.32 of~\cite{cg}.

\begin{lem}\label{univ_cent}
a) An element $x\in \g$ is regular iff $m_i(x)$ form a basis in the centralizer of $x$. 

b) Fix a maximal torus $\ft\subset \g$, then $m_i$ sends  $\ft$ to $\ft$.
Thus we get a regular function $\delta$ on $\ft$ with values in the $\Lambda^r \ft$, $\delta(x) = m_1(x)\wedge \dots \wedge m_r(x)$.
This function has a simple zero on every root hyperplane and no other zeroes.
\end{lem}

In view of Lemma  \ref{univ_cent}a) we get regular functions $\psi_i$ on $\fZ_\g^\g$ such
that for $(x,y)\in \fZ_\g^\g$ ($x\in \g^{reg}$, $y\in \g$, $[x,y]=0$) we have $y=\sum \psi_i(x,y) m_i(x)$. It is clear that function $\phi_i\colon(x,y)\mapsto P_i(x)$
and $\psi_i$, $(i=1,..,r)$ form a coordinate system on $\fZ_\g^\g$.

It is also clear that $\phi_i$ define regular functions on $\fB_\g^\g$, it remains to show the same for $\psi_i$. 
We can consider $\psi_i$ as a rational function on $(\ft\times \ft)/W$, we use the same notation for its pull back to $\ft \times \ft$.
By Cramer rule, Lemma \ref{univ_cent}b) implies that $\psi_i\Delta$ is a regular function on $\ft\times \ft$. Since  $ \psi_i\Delta$ is  $W$-anti-invariant,
it lies in the ideal $I$, thus $\psi_i=\frac{\Delta\psi}{\Delta}$ is a regular function on $\fB_\g^\g$. This proves the part of b) concerned with $\fB_\g^\g$ and $\fZ_\g^\g$.

We now prove statement a). The morphism $\fB_\g^\g\to \ft/W$ is smooth since $\fB_\g^\g\cong \fZ_\g^\g$ as shown above. 
The natural map  $\fBt_\g^\g\to \fB_\g^\g\times_{\ft/W}\ft$ is an isomorphism, since it is finite, birational and its target is smooth (in particular, normal).
 Thus the map $\varpit: \fBt_\g^\g\to\ft$ is smooth.
 We proceed to prove it in the remaining cases by reducing it to the above one.

The morphism $\varpit$ is a composition of an open embedding
and a blow up; let $D$ denote the center of the blow up.
Fix a point $(x,y)$ in one of the spaces $\ft\times T$,  $T\times \ft$, $T\times T$ or $\LT\times T$ (the source of $\varpit$).
Consider the subset $R'$ of roots of $\g$ vanishing on
both $x$ and $y$ (in the case of $\LT\times T$ we use the canonical bijection between the roots of $\g$ and $\Lg$). Let $G'$ be the reductive group with
a maximal torus $T$ and root system $R'$. (Notice that $G'$ is not necessarily identified with a subgroup in $G$, nor is its dual realized as a subgroup in $\LG$;
it is something one can call a bi-endoscopic subgroup of $G$.)
There exists an \'etale neighborhood $U$ of $(x,y)$ such that the pair $(U,D_U)$ is isomorphic to
$(V,D'_{\ft,\ft})$, a neighborhood of 0 in $\ft^2$; here
$D'_\ft$ is the union of root hyperplanes for $\g'$, and $D'_{\ft,\ft}$ is defined the same way as $D_{\ft,\ft}$;
finally $D_U$ denotes the fiber product of $U$ and $D$. Moreover, we can and will assume that
the isomorphism sends the pull back of the discriminant divisor from the second factor (the base) to a divisor
containing $\ft \times D_\ft\cap V$. This shows that  
$\varpit$ is locally isomorphic to the projection $\fBt_{\g'}^{\g'}\to \ft$, in particular it is smooth. 
Thus $\fBt$ is smooth, hence $\fB$ is Cohen-Macaulay in all cases. Since the map $\ft\to \ft/W$, $T\to T/W$ is finite, it follows
that the morphism $\varpit$ is equidimensional. If $G$ is simply-connected, then $T/W$ is also smooth, which implies flatness of $\varpi$.

Now validity of statement  b) in all cases follows from Proposition \ref{blow}, since both sides of the purported isomorphisms are flat over a smooth base
and isomorphic away from a codimension two subvariety of the base.

Furthermore, part d) of the proposition follows as the minimal level for $G$ (viewed as a
$W$-equivariant homomorphism $T\to\LT$) identifies $\LT$ with $T/Z(G)$. Hence
$T\times\LT$ is identified with the quotient of $T\times T$ modulo the action of $Z(G)$
on the second factor. Moreover, this identification takes $D_{T,\LT}$ to
$D_{T,T}/Z(G)$ and is $W$-equivariant.

Similarly, part c) of the proposition follows as the minimal level for $\LG$ (viewed as a
$W$-equivariant homomorphism $\LT\to T$) identifies $T$ with $\LT/Z(\LG)$. Hence
$T\times\LT$ is identified with the quotient of $\LT\times\LT$ modulo the action of
$Z(\LG)$ on the first factor. Moreover, this identification takes $D_{T,\LT}$ to
$Z(\LG)\backslash D_{\LT,\LT}$ and is $W$-equivariant.

Finally, part e) of the proposition follows since the Poisson structures in question
are identified on the open subvariety equal to the preimage of the regular part of the
base (e.g.\ $(T\times \ft^{reg})/W$).

\label{discr}

\section{$W$-invariant sections and blow-ups}

The aim of this section is a proof of Proposition ~\ref{sechenija}. 
We concentrate on the last statement, the other being completely similar.

Let $T^{reg}\subset T$, $T^{reg}_\alpha\subset T$ be the open subschemes defined
by $T^{reg}=\{t\ |\ \alpha(t)\ne 1$ for all roots $\alpha\}$;
  $T^{reg}_\alpha=\{t\ |\ \beta(t)\ne 1$ for all roots $\beta\ne \alpha\}$;
and $T^\bullet=\cupl_\alpha T^{reg}_\alpha$ (thus $T-T^\bullet$ has codimension 2 in $T$
(where the empty subscheme in a curve is considered to be of codimension 2)).
Notice that since $G$ is simply connected the action of $W$ on $T^{reg}$ is free.

We start with a
\begin{lem}
\label{same_hren_lem}
The map $\fBt_G^\LG\times_TT^\bullet \to\fBt_G^\LG/W\times_{T/W}T^\bullet $
is an isomorphism.
\end{lem}
\begin{proof} Let $X\to Y$ be a flat morphism of semi-separated (which means that
the diagonal embedding is affine)
schemes of finite type over a 
characteristic zero field, and let a finite group $W$ act on $X,Y$
so that the map
is $W$-equivariant. Assume that $Y$ is flat over $Y/W$.
We then claim that the  map $X \to X/W \times _{Y/W} Y$ is an 
isomorphism provided that for every Zariski point $y\in Y$ the action of
$\on{Stab}_W(y)$
on the scheme-theoretic fiber $X_y$ is trivial (here $X/W$, $Y/W$ stand for
categorical quotients).
To check this claim we can assume $X$ is affine:
by semi-separatedness every $W$-invariant
subset in $X$ has a $W$-invariant affine neighborhood. Let us
first  assume also that $Y/W$ is a point;
then (by replacing $Y$ by its connected component, and
$W$ by the stabilizer of that component)
 we can assume that $Y$ is nilpotent.
Then $\O_X$ is free over $\O_Y$, and the generators of $\O_X$ as an
$\O_Y$ module can be chosen to be $W$-invariant
(by semi-simplicity of the $W$ action
on $\O_X$, and triviality of the $W$-action on $\O_X\otimes_{\O_Y}{\mathsf k}$);
since $\O_Y^W={\mathsf k}$
(where ${\mathsf k}$ is the base field) we see that
 $\O_X^W \otimes \O_Y\iso\O_X$ as claimed.
Now for a general $Y$ we see that the morphism in question is a morphism
of flat schemes of finite type over $Y/W$,
which induces  an isomorphism on every fiber;
and such a morphism is necessarily an isomorphism.

Now it remains to check that the above conditions hold for
$X= \fBt_G^\LG\times _T T^\bullet $, $Y= T$. For $y\in T^{reg}$ the stabilizer of $y$
is trivial, so there is nothing to check. Consider now $y\in T^{reg}_\alpha $,
$y\not \in T^{reg}$.
Then the stabilizer of $y$ is $\{ 1, s_\alpha\}$.
The ring of functions on $\fBt_G^\LG$ is generated by
$^1\lambdach$, $^2\mu$, $t_\alpha$
 where $\lambdach$, $\mu$ run over weights of $\LT$, $T$ respectively,
$\alpha\in R^+$, and $t_\alpha(^2\alpha-1)=\ ^1\alphach-1$.
We have $s_\alpha^*(^1\lambdach)=\
^1\lambdach\cdot (^1\alphach)^{\langle-\alpha,\lambdach\rangle},\
s_\alpha^*(^2\mu)=\
^2\mu\cdot(^2\alpha)^{\langle-\mu,\alphach\rangle}$,
and $s_\alpha^*(t_\alpha)=t_\alpha\cdot\frac{^2\alpha}{^1\alphach}$.
On the fiber we have $^2\alpha=1$, hence $^1\alphach=1$,
so the action of $s_\alpha$  on the fiber is trivial. \hfill $\Box$
\end{proof}

\bigskip

 Proposition \ref{sechenija} clearly follows from the $(ii) \iff (iv)$ part of
the next

\begin{prop}\label{same_hren}
 Let $S\to T/W$ be a flat morphism, and set $\phi:S\times_{T/W}T^{reg}/W\to
(\LT\times T)/W$ be a $T^{reg}/W$-morphism. Then
 the following are equivalent:

(i) $\phi$ extends to a morphism $S\times_{T/W}T^\bullet /W\to \fB_G^\LG\times_{T/W}T^\bullet $.

(ii)   $\phi$ extends to a morphism $S\to  \fB_G^\LG$.
 
(iii) For every $\alpha\in R$  the morphism
$\phi\times\id_{T^{reg}}: S\times_{T/W} T^{reg}\to \LT \times T^{reg}$
extends to a morphism $S\times_{T/W} T^{reg}_\alpha\to\LT\times T^{reg}_\alpha$
such that ~\eqref{Denis3} holds.

(iv) $\phi\times\id_{T^{reg}}: S\times_{T/W} T^{reg}\to \LT \times T^{reg}$
 extends to a morphism $  S\times_{T/W} T \to \LT\times T $, such that
\eqref{Denis3} holds for every $\alpha\in R$.

\end{prop}

\begin{proof} 
It is enough to assume that $S$ is affine.
Indeed, a morphism from $S$ extends iff
its restriction to every affine open in $S$ does,
because compatibility on intersections
follows from uniqueness of such an extension;
this uniqueness follows from flatness:
if $S$ is flat affine, then tensoring the injection
$\O \to j_*\O$ with $\O_S$ we get
an imbedding $\O_S\hookrightarrow j_*j^*\O_S$,
where $j$ stands for the imbedding
 $T^{reg}/W\to T/W$, or $T^{reg}\to T$.
So we will assume $S$ affine from now on.

(iv) $\Rightarrow$ (iii) and (ii) $\Rightarrow$ (i) are obvious.

 To check that (iii)  $\Rightarrow$ (iv) 
we tensor (over $\O_{T/W}$) the exact sequence of   $\O_T$-modules
\begin{equation}\label{sequence}
0\to \O \to  \O_{T^{reg}} \to \oplusl_\alpha (\O_{T^{reg}_\alpha}/\O_T)
\end{equation}
with $\O_S$. The resulting exact sequence shows that  a regular function on $S\times _{T/W}
T^{reg}$ extends to a regular function on $S\times_{T/W}T$ iff it extends to $S\times _T
T^{reg}_\alpha$ for all $\alpha$. Applying this observation to
 $(\phi\times\id)^*(f|_{ \LT\times T^{reg}})$
 for each regular
function $f$ on $ \LT\times T $ we see that ~(iii) implies extendability of
$\phi\times\id$ to $ S\times_{T/W} T$.
It is also clear that ~\eqref{Denis3} holds if it holds
on $T^\bullet$.

Verification of (i)  $\Rightarrow$ (ii) is similar (with ~\eqref{sequence}
replaced by the $W$-invariant part of ~\eqref{sequence}).


It remains to check (i) $\iff$ (iii).
If (i) holds, i.e.
$\phi$ extends to a map $ S\times_{T/W}T^\bullet /W\to \fB_G^\LG\times_{T/W}T^\bullet $
then we can take the fiber product of this map with $\id_{T^\bullet}$ over $T/W$.
By Lemma ~\ref{same_hren_lem} it yields a map
$S\times_{T/W}T^\bullet\to\fBt_G^\LG\times_TT^\bullet$,
which can be composed with the projection $\fBt_G^\LG\to\LT\times T$
to produce a map
$ S\times_{T/W}T^\bullet\to\LT\times T^\bullet$.
It is clear that this map satisfies ~\eqref{Denis3},
because the image of the map $\fBt_G^\LG\to  \LT\times T$
intersected with $\LT\times\on{Ker}(^2\alpha)$
is contained in $\on{Ker}(^1\alphach)\times T$.

Conversely, if (iii) holds then restricting the given map $S\times _{T/W}T^\bullet
\to \LT\times T^\bullet$ to $S\times _{T/W} (\on{Ker}(\alpha)\cap T^\bullet)$ we get
a map into $\on{Ker}(\alphach)\times T$
(this is immediate from ~\eqref{Denis3}).
This means that the map lifts to a map into $\fBt_G^\LG$. Replacing both
the source and the target by their  quotients
by $W$ we get the map required in ~(i). \hfill $\Box$
\end{proof}

\section{$K$-theory and blow-ups}
The aim of this section is a proof of Proposition ~\ref{K}. Recall that
~\ref{K} ~(a) was already proved in ~\ref{kommut}. $G$ is assumed
simply connected till the further notice.

\subsection{Reminder on the affine Grassmannians}
\label{reminder}
Let $X=X_G$ be the lattice of characters of $T$, and let
$Y=Y_G$ be the lattice of cocharacters of $G$. Note that $X_G=Y_\LG,\
Y_G=X_\LG$. Let $X^+\subset X$ (resp. $Y^+\subset Y$) be the cone of
dominant weights (resp. dominant coweights). It is well known that
the $\GO$-orbits in $\Gr_G$ are numbered by the dominant coweights:
$\Gr_G=\bigsqcup_{\lambdach\in Y^+}\Gr_{G,\lambdach}$. The adjacency relation
of orbits corresponds to the standard partial order on coweights:
$\ol\Gr_{G,\lambdach}=\bigsqcup_{\much\leq\lambdach}\Gr_{G,\much}$.
The open embedding $\Gr_{G,\lambdach}\hookrightarrow\ol\Gr_{G,\lambdach}$
will be denoted by $j_\lambdach$ or simply by $j$ if no confusion is likely.
The dimension $\dim(\Gr_{G,\lambdach})=\langle2\rho,\lambdach\rangle$
where $2\rho=\sum_{\alpha\in R^+}\alpha$, and
$\langle,\rangle:\ X\times Y\to\BZ$ is the canonical perfect pairing.

Recall that the $T$-fixed points in $\Gr_G$ are naturally numbered by
$Y$; a point $\much$ lies in an orbit $\Gr_{G,\lambdach}$ iff $\much$ lies
in the $W$-orbit of $\lambdach$. Each $\GO$-orbit $\Gr_{G,\lambdach}$ is
partitioned into Iwahori orbits isomorphic to affine spaces and numbered
by $\much\in W\lambdach$. Hence the basics of ~\cite{cg} ~Chapter ~5 are
applicable in our situation.

In particular, $K^T(\Gr_{G,\lambdach})$ is a free $K^T(pt)$-module, and
$K^\GO(\Gr_{G,\lambdach})=K^G(\Gr_{G,\lambdach})$
is a free $K^G(pt)$-module
(recall that $K^T(pt)=\BC[T]$, and $K^G(pt)=\BC[T/W]$).
Moreover, the natural map $K^T(pt)\otimes_{K^G(pt)}K^G(\Gr_{G,\lambdach})\to
K^T(\Gr_{G,\lambdach})$ is an isomorphism, and
$K^G(\Gr_{G,\lambdach})=K^T(\Gr_{G,\lambdach})^W$, cf. ~\cite{cg} ~6.1.22.

Since $K^\TO(\Gr_G)=K^T(\Gr_G)$ (resp. $K^\GO(\Gr_G)=K^G(\Gr_G)$) is filtered
by the support in $\GO$-orbit closures, with the associated graded
$\bigoplus_{\lambdach\in Y^+}K^T(\Gr_{G,\lambdach})$ (resp.
$\bigoplus_{\lambdach\in Y^+}K^G(\Gr_{G,\lambdach})$), we arrive at the
following

\begin{lem}
\label{flat 2}
$K^\TO(\Gr_G)=K^T(\Gr_G)$ is a flat $K^T(pt)$-module, and
$K^\GO(\Gr_G)=K^G(\Gr_G)$ is a flat $K^G(pt)$-module.
Moreover, the natural map $K^T(pt)\otimes_{K^G(pt)}K^G(\Gr_G)\to
K^T(\Gr_G)$ is an isomorphism, and $K^G(\Gr_G)=\left(K^T(\Gr_G)\right)^W$.
\end{lem}

\subsection{Localization}
\label{loc}
The space $K^T(\Gr_G)=K^\TO(\Gr_G)=K(\TO\backslash\GF/\GO)$ is equipped with
the two commuting actions: $K(\TO\backslash\TF/\TO)$ acts by convolutions
on the left, and $K^G(\Gr_G)=K^\GO(\Gr_G)=K(\GO\backslash\GF/\GO)$
acts by convolutions on the right. Also, $W$ acts on $K^T(\Gr_G)$ commuting
with the right action of $K^G(\Gr_G)$.
Clearly, the algebra $K(\TO\backslash\TF/\TO)$ is isomorphic to
$\BC[\LT\times T]$. The action of $W$ on $K^T(\Gr_G)$ normalizes the action of
$K(\TO\backslash\TF/\TO)$ and induces the natural (diagonal) action of $W$
on $\BC[\LT\times T]$.

Let $g$ be a general (regular) element of $T$. Then the fixed point set
$(\Gr_G)^g=(\Gr_G)^T=Y$ coincides with the image of the embedding
$\Gr_T\hookrightarrow\Gr_G$. According to Thomason Localization Theorem
(see e.g. ~\cite{cg} ~5.10), after localization,
$\left(K^T(\Gr_G)\right)_g$
becomes a free rank one $\left(K(\TO\backslash\TF/\TO)\right)_g$-module.
This means that after restriction to $T^{reg}\subset T=\on{Spec}(K^T(pt))$
we have an isomorphism $K^T(\Gr_G)|_{T^{reg}}\simeq\BC[\LT\times T]|_{T^{reg}}$
compatible with the natural $W$-actions. The localized algebra
$K^\GO(\Gr_G)|_{T^{reg}/W}$ is embedded into
$\left(\End_{K(\TO\backslash\TF/\TO)|_{T^{reg}}}
(K^T(\Gr_G)|_{T^{reg}})\right)^W$.
According to Lemma ~\ref{flat 2}, $K^G(\Gr_G)=\left(K^T(\Gr_G)\right)^W$;
hence this embedding is an isomorphism, and we have
$K^\GO(\Gr_G)|_{T^{reg}/W}\simeq\BC[\LT\times T]^W|_{T^{reg}/W}$.

Hence both $\BC[\fB_G^\LG]$ and
$K^\GO(\Gr_G)$ are the flat $\BC[T]^W$-modules embedded into
$\BC[\LT\times T](\Delta^{-1})$ (see ~\ref{discr}).
We must prove that the identification of
$\BC[\fB_G^\LG]|_{T^{reg}/W}$
and $K^\GO(\Gr_G)|_{T^{reg}/W}$ extends to the identification
over the whole $T/W$. To this end it suffices to check that the identification
extends over the codimension 1 points of $T/W$. Let $g\in T/W$ be a regular
point of $\bD$; that is, $g$ is represented by a semisimple element of $G$
such that the centralizer $Z(g)$ has semisimple rank 1.

We must prove that the localizations $\BC[\fB_G^\LG]_g$ and
$\left(K^\GO(\Gr_G)\right)_g$
are isomorphic. To this end it suffices to identify
$\BC\left[\LT\times T,\
\frac{^1\alphach-1}{^2\alpha-1},\ \alpha\in R\right]_g$
(which we denote by $\BC[\overset{\bullet}\fB{}_G^\LG]_g$ for short)
and $\left(K^T(\Gr_G)\right)_g$.
Note that the embedding of reductive groups
$Z(g)\hookrightarrow G$
(the neutral connected component) induces the isomorphism
$\Gr_{Z(g)}=(\Gr_G)^g\hookrightarrow\Gr_G$. According to Thomason Localization
Theorem, we have an isomorphism of localizations
$\left(K^T(\Gr_{Z(g)})\right)_g\simeq\left(K^T(\Gr_G)\right)_g$.
Finally, the isomorphism
$K^T\left(\Gr_{Z(g)}\right)\simeq
\BC\left[\overset{\bullet}\fB{}_{Z(g)}^{\check{Z}(g)}\right]$
follows from the
calculations in ~\ref{sootn}, ~\ref{Iwahori},
and together with the evident isomorphism of localizations
$\BC\left[\overset{\bullet}\fB{}_{Z(g)}^{\check{Z}(g)}\right]_g\simeq
\BC\left[\overset{\bullet}\fB{}_G^\LG\right]_g$
establishes the desired isomorphism
$\left(K^T(\Gr_G)\right)_g\simeq\BC\left[\overset{\bullet}\fB{}_G^\LG\right]_g$.

This completes the proof of ~\ref{K} ~(b).

\subsection{Comparison of Poisson structures}
In order to compare the Poisson structures on $K^\GO(\Gr_G)$
and $\BC[\fB_G^\LG]$ it suffices to identify them on the open subset
$K^\GO(\Gr_G)|_{T^{reg}/W}=\BC[\fB_G^\LG]|_{T^{reg}/W}=
\BC[\LT\times T^{reg}]^W$.
The space $$K^{T\times\Gm}(\Gr_G)=K^{\TO\ltimes\Gm}(\Gr_G)=
K\left(\TO\ltimes\Gm\backslash\GF\ltimes\Gm/\GO\ltimes\Gm\right)$$
is equipped with
the two commuting actions:
$K\left(\TO\ltimes\Gm\backslash\TF\ltimes\Gm/\TO\ltimes\Gm\right)$
acts by convolutions
on the left, and $$K^{G\times\Gm}(\Gr_G)=K^{\GO\ltimes\Gm}(\Gr_G)=
K\left(\GO\ltimes\Gm\backslash\GF\ltimes\Gm/\GO\ltimes\Gm\right)$$
acts by convolutions on the right.
Also, $W$ acts on $K^{\TO\ltimes\Gm}(\Gr_G)$ commuting
with the right action of $K^{\GO\ltimes\Gm}(\Gr_G)$.
Clearly, the algebra
$K\left(\TO\ltimes\Gm\backslash\TF\ltimes\Gm/\TO\ltimes\Gm\right)$
is isomorphic to the group algebra $\BC[\Gamma]$ of the following
Heisenberg group $\Gamma$.

It is a $\BZ$-central extension of $Y\times X$ with the multiplication
(written multiplicatively)
$$(q^{n_1},e^{\lambdach_1},e^{\mu_1})\cdot(q^{n_2},e^{\lambdach_2},e^{\mu_2})=
(q^{n_1+n_2+\langle\mu_1,\lambdach_2\rangle},e^{\lambdach_1+\lambdach_2},
e^{\mu_1+\mu_2})$$
where $\langle,\rangle:\ X\times Y\to\BZ$ is the canonical perfect pairing.

Finally, the action of the Weyl group $W$ on $K^{\TO\ltimes\Gm}(\Gr_G)$
normalizes the action of
$K\left(\TO\ltimes\Gm\backslash\TF\ltimes\Gm/\TO\ltimes\Gm\right)$
and induces the natural (diagonal) action of $W$
on $\BC[\Gamma]$. From this we deduce, exactly as in ~\ref{loc}, that
$K^{\GO\ltimes\Gm}(\Gr_G)|_{T^{reg}/W}\simeq\BC[\Gamma]|_{T^{reg}/W}$.
It follows that the Poisson structure on $K^\GO(\Gr_G)|_{T^{reg}/W}$
coincides with the standard Poisson structure on $\BC[\LT\times T^{reg}]^W$.

This completes the proof of ~\ref{K} ~(c).

\subsection{The case of non simply connected $G$}
For general $G$ let $\tilde G$ denote its universal cover,
and let $\tilde T$ stand for the Cartan of $\tilde G$.
Note that the dual torus is $\LT/\pi_1(G)$.
As in ~\ref{loc}, we have
$K^G(\Gr_G)=\left(\End_{K(\TO\backslash\TF/\TO)}(K^T(\Gr_G))\right)^W$,
so it suffices to identify the
$K(\TO\backslash\TF/\TO)\ltimes W=\BC[\LT\times T]\ltimes W$-module
$K^T(\Gr_G)$ with
$\BC\left[\LT\times T,\
\frac{^1\alphach-1}{^2\alpha-1},\ \alpha\in R\right]
=\on{Spec}\BC[\overset{\bullet}\fB{}_G^\LG]$.
We do this by reduction to the known case of $\tilde G$.

Evidently, the
$K(\TO\backslash\TF/\TO)\ltimes W=\BC[\LT\times T]\ltimes W$-module
$K^T(\Gr_G)$ equals
$\BC[\LT\times T]\ltimes W\otimes_{\BC[(\LT/\pi_1(G))\times T]\ltimes W}
K^T(\Gr_{\tilde G})$. On the other hand, it follows from ~\ref{loc} that the
$K(\TO\backslash{\tilde T}(\bF)/{\tilde T}(\bO))\ltimes W=
\BC[(\LT/\pi_1(G))\times T]\ltimes W$-module
$K^T(\Gr_{\tilde G})$ equals
the invariants of $\pi_1(G)$ in $K^{\tilde T}(\Gr_{\tilde G})$, that is
$\BC\left[(\LT/\pi_1(G))\times{\tilde T},\
\frac{^1\alphach-1}{^2\alpha-1},\ \alpha\in R\right]^{\pi_1(G)}=
\BC\left[(\LT/\pi_1(G))\times T,\
\frac{^1\alphach-1}{^2\alpha-1},\ \alpha\in R\right]$.

This completes the proof of ~\ref{K} for general $G$.

\subsection{Borel-Moore Homology and blow-ups}
\label{parallel}
Theorem ~\ref{HBM} is proved absolutely parallelly to the proof of
Theorem ~\ref{K}.

\section{Computation of $K_\GO(\Lambda)$.}
\label{seventeen}

\subsection{The affine Grassmannian Steinberg variety}
\label{steinberg}
We denote by $\fu\subset\fg(\bO)$ (resp. $U\subset\GO$) the nilpotent
(resp. unipotent) radical. It has a filtration 
$\fu=\fu^{(0)}\supset\fu^{(1)}\supset\ldots$ 
by congruence subalgebras. 
The trivial (Tate) vector bundle $\ul{\fg(\bF)}$ with the fiber
$\fg(\bF)$ over $\Gr_G$ has a structure of an ind-scheme. It
contains a profinite dimensional vector subbundle
$\ul\fu$ whose fiber over a point $g\in\Gr_G$ represented by a compact
subalgebra in $\fg(\bF)$ is the pronilpotent radical of this subalgebra.
The trivial vector bundle $\ul{\fg(\bF)}=\fg(\bF)\times\Gr_G$ also contains
a trivial vector subbundle $\fu\times\Gr_G$. 

We will call $\ul\fu$ {\em the cotangent bundle} of $\Gr_G$, and we will
call the intersection $\Lambda:=\ul\fu\cap(\fu\times\Gr_G)$ {\em the affine
Grassmannian Steinberg variety}. It has a structure of an ind-scheme of
ind-infinite type. Namely, if $p$ stands for the natural projection 
$\Lambda\to\Gr_G$, then $\Lambda_{\leq\lambdach}:=p^{-1}(\ol\Gr_{G,\lambdach})$
is a scheme of infinite type, and $\Lambda=\bigcup\Lambda_{\leq\lambdach}$.

Note that for a fixed $\lambdach$ and $l\gg0$ the intersection of fibers
of $\ul\fu$ over all points of $\ol\Gr_{G,\lambdach}$ (as vector subspaces
of $\fg(\bF)$) contains $\fu^{(l)}$. Thus $\fu^{(l)}$ acts freely (by fiberwise
translations) on $\Lambda_{\leq\lambdach}$, and the quotient is a scheme of
finite type, to be denoted by $\Lambda_{\leq\lambdach}^l$. 
For $k>l$ we have evident
affine fibrations 
$p^k_l:\ \Lambda_{\leq\lambdach}^k\to\Lambda_{\leq\lambdach}^l$,
and $\Lambda_{\leq\lambdach}$ coincides with the inverse limit of this system.

Similarly, the total space of the vector bundle $\ul\fu$ (to be denoted by
the same symbol) is a union of infinite type schemes $\ul\fu_{\leq\lambdach}$,
and for fixed $\lambdach$ and $l\gg0$, the scheme $\ul\fu_{\leq\lambdach}$ is
the inverse limit of affine fibrations $p^k_l:\ \ul\fu_{\leq\lambdach}^k\to
\ul\fu_{\leq\lambdach}^l\ (k>l)$. Note that the proalgebraic group $\GO$ acts
on all the above schemes, and the fibrations $p^k_l$ are $\GO$-equivariant.

A {\em $\GO$-equivariant coherent sheaf $\CF$ on $\ul\fu$} is by definition
supported on some $\ul\fu_{\leq\lambdach}$. 
There, it is defined as a collection
of $\GO$-equivariant sheaves $\CF^l$ on $\ul\fu_{\leq\lambdach}^l$ for $l\gg0$
together with isomorphisms $(p^k_l)^*\CF^l\simeq\CF^k$. We will consider
the $\GO$-equivariant coherent sheaves on $\ul\fu$ {\em supported on 
$\Lambda$}, and $D^bCoh^\GO_\Lambda(\ul\fu)$ 
stands for the derived category of 
such sheaves, and $K^\GO(\Lambda)$ stands for the $K$-group of such sheaves.

\subsection{Convolution in $D^bCoh^\GO_\Lambda(\ul\fu)$}
\label{convolu}
We have a principal $\GO$-bundle $\GF\to\Gr_G$. Given a $\GO$-(ind)-scheme
$A$ we can form an associated bundle $\widetilde{A}=\GF\times_\GO A\to\Gr_G$. 
Given a coherent $\GO$-equivariant sheaf $\CF$ on $A$ we can form an
associated sheaf $\widetilde\CF$ on $\widetilde{A}$ as $\GO$-invariants
in the direct image of $\CO_{\GF}\boxtimes\CF$ from $\GF\times A$ to
$\GF\times_\GO A$. If $A=\Gr_G$, apart from the natural projection
$p_1:\ \widetilde{A}\to\Gr_G$, we have a multiplication map
$\GF\times_\GO\Gr_G\to\Gr_G$, to be denoted $p_2$. Then $(p_1,p_2)$ identifies
$\widetilde{\Gr_G}$ with $\Gr_G\times\Gr_G$. Furthermore, 
$\widetilde{\ul\fu}$ is a vector bundle over 
$\widetilde{\Gr_G}=\Gr_G\times\Gr_G$ which is naturally identified with
$p_2^*\ul\fu$. Thus we have an ind-proper morphism 
$p_2:\ \widetilde{\ul\fu}\to\ul\fu$.

Note that both $\widetilde{\ul\fu}=p_2^*\ul\fu$ and $p_1^*\ul\fu$ are
subbundles in the trivial (Tate) vector bundle $\ul{\fg(\bF)}$ 
over $\Gr_G\times\Gr_G$ with the fiber $\fg(\bF)$. Their intersection
is naturally identified with $\widetilde\Lambda$. In particular, we have
an embedding $\widetilde{\Lambda}\subset p_1^*\ul\fu\oplus p_2^*\ul\fu$,
and an ind-proper morphism $p_2:\ \widetilde{\Lambda}\to\ul\fu$.

Hence given $\GO$-equivariant coherent sheaves $\CF,\CG$ on $\Lambda$
we can consider the $\GO$-equivariant complex $\CF\star\CG:=
(p_2)_*(p_1^*\CF\stackrel{L}{\otimes}\widetilde\CG)$ (tensor product over the
structure sheaf of the profinite dimensional vector bundle
$p_1^*\ul\fu\oplus p_2^*\ul\fu$). Clearly, 
$\CF\star\CG$ is supported on $\Lambda$. Hence we get
a convolution operation on $D^bCoh^\GO_\Lambda(\ul\fu)$ and on $K^\GO(\Lambda)$
once we check that $p_1^*\CF\stackrel{L}{\otimes}\widetilde\CG$
is bounded. 

To this end, note that $\widetilde\CG$ is flat over the first copy of
$\Gr_G$, and for some $\lambdach$ the sheaf $\CF$ is supported on
$\Lambda_{\leq\lambdach}$, so the tensor product
$p_1^*\CF\stackrel{L}{\otimes}\widetilde\CG$ can actually be computed
over the structure sheaf of 
$(p_1^*\ul\fu\oplus p_2^*\ul\fu)|_{\ol\Gr_{G,\lambdach}\times\Gr_G}=
\ul\fu_{\leq\lambdach}\times\ul\fu\subset\ul\fu\times\ul\fu=
p_1^*\ul\fu\oplus p_2^*\ul\fu$. 
That is, $p_1^*\CF\stackrel{L}{\otimes}\widetilde\CG$ is the direct image of
$p_1^*\CF|_{\ul\fu_{\leq\lambdach}\times\ul\fu}
\stackrel{L}{\otimes}_{\CO_{\ul\fu_{\lambdach}\times\ul\fu}}
\widetilde\CG|_{\ul\fu_{\leq\lambdach}\times\ul\fu}$ under the closed
embedding $\ul\fu_{\leq\lambdach}\times\ul\fu\hookrightarrow
\ul\fu\times\ul\fu$.
On the other hand, $p_1^*\CF$ is flat over the second copy of $\Gr_G$,
while the support of $\widetilde\CG$ intersected with 
$\ul\fu_{\leq\lambdach}\times\ul\fu$
is contained in
$\ul\fu_{\leq\lambdach}\times\ul\fu_{\leq\much}$
for some $\much$. Hence the tensor product
$p_1^*\CF\stackrel{L}{\otimes}\widetilde\CG$ can actually be computed
over the structure sheaf of 
$\ul\fu_{\leq\lambdach}\times\ul\fu_{\leq\much}$.
There exists $l\gg0$ such that the diagonal fiberwise action of 
$\fu^{(l)}$ on $\ul\fu_{\leq\lambdach}\times\ul\fu_{\leq\much}$ is free,
and both $p_1^*\CF$ and $\widetilde\CG$ restricted to 
$\ul\fu_{\leq\lambdach}\times\ul\fu_{\leq\much}$ are $\fu^{(l)}$-equivariant,
that is, they are lifted from the sheaves on 
$(\ul\fu_{\leq\lambdach}\times\ul\fu_{\leq\much})/\fu^{(l)}=:V$;
we abuse notation by keeping the same names for these sheaves.
So the tensor product
$p_1^*\CF\stackrel{L}{\otimes}\widetilde\CG$ can actually be computed
as the tensor product of coherent sheaves
over the structure sheaf of the profinite dimensional vector bundle $V$
over the finite dimensional scheme $\ol\Gr_{G,\lambdach}\times
\ol\Gr_{G,\much}$.

Now there exists a vector subbundle $V'\subset V$ 
such that the quotient $\ol{V}:=V/V'$ is a finite dimensional
vector bundle, $p_1^*\CF$ is lifted from $\ol V$, 
and the support of $\widetilde\CG$ in $V$ projects isomorphically
onto its image in $\ol V$. Moreover, recall that $p_1^*\CF$ is flat
over $\ol\Gr_{G,\much}$, while $\widetilde\CG$ is flat over 
$\ol\Gr_{G,\lambdach}$. Clearly, in this situation
$p_1^*\CF\stackrel{L}{\otimes}\widetilde\CG\in D^b(V)$.
This explains why for $\GO$-equivariant coherent sheaves $\CF,\CG$
on $\Lambda$ the tensor product
$p_1^*\CF\stackrel{L}{\otimes}\widetilde\CG$ is a bounded complex of
coherent sheaves on $p_1^*\ul\fu\oplus p_2^*\ul\fu$ supported on
$\widetilde\Lambda$. Hence the same is true for the bounded complexes
of $\GO$-equivariant coherent sheaves $\CF,\CG$ on $\ul\fu$ supported
on $\Lambda$. Thus, $D^bCoh^\GO_\Lambda(\ul\fu)$ is closed with respect to 
convolution.

\begin{thm}
\label{none}
$K^\GO(\Lambda)$ is a commutative algebra isomorphic to 
$\BC[\check{T}\times T]^W$.
\end{thm}

\begin{rem}\label{cherednik}
{\em Since $\Lambda_G$
is an affine Grassmannian analogue of the classical Steinberg variety,
this result agrees well with the geometric realization of the Cherednik
double affine Hecke algebra in ~\cite{gg}, ~\cite{v}. In effect, 
$K^\GO(\Lambda_G)$ is the spherical subalgebra of the Cherednik algebra
with both parameters trivial: $q=t=1$.}
\end{rem}

\subsection{Bialynicki-Birula stratifications}
\label{bb}
The proof of Theorem ~\ref{none} uses the following lemma on 
$K$-theory of cellular spaces. Let $M$ be a normal quasiprojective variety
equipped with a torus $H$-action with finitely many fixed points. 
We assume that $M$ is equipped with an $H$-invariant stratification
$M=\bigsqcup_{\mu\in M^H}M_\mu$ such that each stratum $M_\mu$ contains
exactly one $H$-fixed point $\mu$, and $M_\mu$ is isomorphic to an affine
space.
For $\mu\in M^H$ we denote by $j_\mu:\ M_\mu\hookrightarrow M$ the locally
closed embedding of the corresponding stratum.
We denote by $i_\mu:\ \mu\hookrightarrow M_\mu$ the closed embedding of an 
$H$-fixed point in the corresponding stratum, or in the whole of $M$ when no
confusion is likely. We denote by $\mu\leq\nu$ the
closure relation of strata. We denote by $M_{\leq\mu}\subset M$ the 
union $\bigcup_{\nu\leq\mu}M_\nu$.

Given an $H$-equivariant closed embedding of $M$ into a smooth
$H$-variety $M'$ (for the existence see ~\cite{su})
we denote by $T^*M$ the restriction of the cotangent bundle
$T^*M'$ to $M\subset M'$.
We denote by $\imath:\ M\hookrightarrow T^*M$ the embedding of the zero
section. 
We also denote by $i_\mu$ the closed embedding of the conormal bundle
$T^*_\mu M'\hookrightarrow T^*M$ when no confusion is likely.
Finally, we denote by $\CL'$ the union of conormal bundles
$\bigcup_\mu T^*_{M_\mu}M'$, and $\jmath$ stands for the closed embedding
$\CL'\hookrightarrow T^*M$. We denote by $\CL'_{\leq\mu}\subset\CL'$ the
union $\bigcup_{\nu\leq\mu}T^*_{\M_\nu}M'$; it is a closed subvariety
of $\CL'$. It has a closed subvariety 
$\CL'_{<\mu}:=\bigcup_{\nu<\mu}T^*_{M_\nu}M'$.

For $\mu\in M^H$ we have an embedding
$i_{\mu*}:\ K^H(\mu)\hookrightarrow K^H(M)$. 
We have an embedding $\jmath_*:\ K^H(\CL')\hookrightarrow K^H(T^*M)
\stackrel{\imath^*}{\simeq}K^H(M)$. 
Indeed, the exact sequences (see ~\cite{cg} Chapter 5)
$$0\to K^H(\CL'_{<\mu})\to K^H(\CL'_{\leq\mu})\to K^H(T^*_{M_\mu}M')\to0,$$
$$0\to K^H(T^*M'|_{M_{<\mu}})\to K^H(T^*M'|_{M_{\leq\mu}})\to 
K^H(T^*M'|_{M_\mu})$$
give rise to the support filtrations on $K^H(\CL')$ and $K^H(T^*M)$
with associated graded $\bigoplus_{\mu\in M^H}K^H(T^*_{M_\mu}M')$
and $\bigoplus_{\mu\in M^H}K^H(T^*M'|_{M_\mu})$. Now $\jmath_*$ is strictly
compatible with the support filtrations and clearly injective on the
associated graded.

Note that the image $\jmath_*(K^H(\CL'))
\subset K^H(M)$ is independent of the choice of the closed embedding
$M\hookrightarrow M'$. In effect, given another embedding 
$M\hookrightarrow\widetilde M$, we can consider the diagonal embedding
$M\hookrightarrow M'':=M'\times\widetilde M$. Clearly, we have a projection
$p:\ T^*M''|_M\to T^*M'|_M$ which realizes $T^*M''|_M$ as a vector bundle
over $T^*M'|_M$. Moreover, if we denote by $\CL''$ the union of
conormal bundles $\bigcup_\mu T^*_{M_\mu}M''\subset T^*M''|_M$ then
$\CL''=p^{-1}\CL'$. This shows that the images of $K^H(\CL')$ and
$K^H(\CL'')$ in $K^H(M)$ coincide, and thus $\jmath_*(K^H(\CL'))\subset K^H(M)$
is well-defined.

\begin{lem}
\label{BB}
In $K^H(M)$ we have an equality 
$\jmath_*(K^H(\CL'))=\oplus_\mu i_{\mu*}(K^H(\mu))$.
\end{lem}

\begin{proof}
Let $K^H(D_M)$ stand for the $K$-group of weakly $H$-equivariant $D$-modules
on $M'$ supported on $M\subset M'$. 
Given such a $D$-module and passing to associated graded with respect
to a good filtration, we obtain an $H$-equivariant coherent sheaf on $T^*M$,
and this way one obtains a homomorphism 
$SS:\ K^H(D_M)\to K^H(T^*M)\stackrel{\imath^*}{\simeq}K^H(M)$
(see e.g. ~\cite{vc}).
Let $\delta_\mu$ stand for a $\delta$-function $D$-module at the point
$\mu\in M^H$ with its obvious $H$-equivariance. Then, evidently,
$SS(\delta_\mu)$ generates $i_{\mu*}(K^H(\mu))$ as a module over
$K^H(pt)$. Moreover, $\{SS(j_{\mu!}\CO_{M_\mu}),\ \mu\in M^H\}$ forms a
basis of $\jmath_*(K^H(\CL'))$. 

In effect, the closed embedding $\CL'_{<\mu}\hookrightarrow\CL'_{\leq\mu}$
gives rise to the exact sequence 
$$0\to K^H(\CL'_{<\mu})\to K^H(\CL'_{\leq\mu})\to K^H(T^*_{M_\mu}M')\to0$$ 
(see ~\cite{cg} Chapter 5), and the image of $SS(j_{\mu!}\CO_{M_\mu})$
in $K^H(T^*_{M_\mu}M')$ clearly generates it.

So it is enough to check the equality in $K^H(T^*M)$:

\begin{equation}
\label{ss}
SS(\delta_\mu)=SS(j_{\mu!}\CO_{M_\mu})\cdot(-1)^{\dim M_\mu}
\det(T_\mu M_\mu)
\end{equation}
where $\det(T_\mu M_\mu)$ is the character of $H$ (thus an 
invertible element of $K^H(pt)=\BC[H]$) acting in the determinant of the
tangent bundle of $M_\mu$ at $\mu$.

To this end note that restriction to the $H$-fixed points gives rise to
an embedding $\oplus_\nu i_\nu^*\imath^*:\ 
K^H(T^*M)\hookrightarrow\oplus_\nu K^H(\nu)$. 
This is checked by induction in $\nu$ using the exact sequences
$$0\to K^H(T^*M'|_{M_{<\nu}})\to K^H(T^*M'|_{M_{\leq\nu}})\to
K^H(T^*M'|_{M_\nu})\to0.$$
It is clear that for $\nu=\mu$
the restrictions $i_\mu^*\imath^*$ of the LHS and RHS of ~(\ref{ss})
coincide. We are going to check that for $\nu\ne\mu$ the restrictions
$i_\nu^*\imath^*$ of the LHS and RHS of ~(\ref{ss}) both vanish.
Evidently, $i_\nu^*\imath^*SS(\delta_\mu)=0$.

Recall that $i_\nu$ also stands for the closed embedding 
$T^*_\nu M'\hookrightarrow T^*M$, so we just have to check that
$i_\nu^*SS(j_{\mu!}\CO_{M_\mu})=0\in K^H(T^*_\nu M')$. Note that the 
functor of global
sections of $H$-equivariant coherent sheaves on the vector space
$T^*_\nu M'$ gives rise to an embedding 
$\Gamma:\ K^H(T^*_\nu M')\hookrightarrow
\BZ^{X^*(H)}$ where $X^*(H)$ stands for the lattice of characters of $H$.
Now for a $D$-module $\CF$ we have $\Gamma(i_\nu^*SS\CF)={\mathbf i}_\nu^*\CF$
where ${\mathbf i}_\nu^*\CF$ stands for the fiber at $\nu\in M$ of the
$H$-equivariant quasicoherent $\CO_{M'}$-module $\CF$. Finally, for 
$\CF=j_{\mu!}\CO_{M_\mu}$ and $\nu\ne\mu$ we have 
${\mathbf i}_\nu^*j_{\mu!}\CO_{M_\mu}=0$. This completes the proof of the
lemma.

\end{proof}

\subsection{Bialynicki-Birula stratification of $\Gr_G$}
\label{BiBi}
We consider the stratification of $\Gr_G$ by the Iwahori orbits
$\Gr_G=\bigsqcup_{\much\in Y}\Gr_G^\much$. This is a refinement of the
stratification by the $\GO$-orbits: 
$\Gr_{G,\lambdach}=\bigsqcup_{\much\in W\lambdach}\Gr_G^\much$. 
Let us denote by $\fn\supset\fu$ the nilpotent radical of the Iwahori
subalgebra in $\fg(\bF)$. The union of conormal bundles to the Iwahori
orbits is the following subvariety $\Lambda_I$ of the cotangent bundle
$\ul\fu$: by definition, $\Lambda_I:=\ul\fu\cap(\fn\times\Gr_G)$.
We have a closed embedding $\Lambda\subset\Lambda_I$.

Lemma ~\ref{BB} allows us to compute 
$K^T(\Lambda_I)=\oplus_{\much\in Y}K^T(\much)\subset K^T(\Gr_G)$,
i.e. $K^T(\Lambda_I)\simeq\BC[\LT\times T]$ (note that the natural
$W$-action on $K^T(\Gr_G)$ induces the diagonal $W$-action on
$\BC[\LT\times T]\simeq K^T(\Lambda_I)\subset K^T(\Gr_G)$). Although
Lemma ~\ref{BB} was formulated for finite dimensional varieties $M$, 
its proof goes through for $\Gr_G$ without changes: we only need to 
have the singular support map $SS: K^T(D_{\Gr_G})\to K^T(\ul\fu)\simeq
K^T(\Gr_G)$. For this see ~\cite{kt},~\cite{bd} (Chapter 15),~\cite{gg}.

The embedding $\Lambda\hookrightarrow\Lambda_I$ gives rise to the
embedding $K^T(\Lambda)\hookrightarrow K^T(\Lambda_I)\hookrightarrow
K^T(\ul\fu)=K^T(\Gr_G)$. Note that $W$ acts naturally on both
$K^T(\Lambda)$ and $K^T(\Gr_G)$, and the embedding
$K^T(\Lambda)\hookrightarrow K^T(\Gr_G)$ is $W$-equivariant.
Also, $(K^T(\Lambda))^W=K^G(\Lambda)=K^\GO(\Lambda)$.
Hence, the image of the embedding $K^\GO(\Lambda)\hookrightarrow
K^T(\Lambda_I)\simeq\BC[\LT\times T]\subset K^T(\Gr_G)$ lies in the
invariants of the diagonal $W$-action on $\BC[\LT\times T]$.
Thus to prove Theorem ~\ref{none} we must check that the image
of this embedding contains $\BC[\LT\times T]^W$.

We have projections $\pi:\ \Lambda\to\Gr_G$, and $\pi_I:\ \Lambda_I\to\Gr_G$.
For $\lambdach\in Y^+$ we denote by $\Lambda_\lambdach$ (resp. 
$\Lambda_{\leq\lambdach},\ \Lambda_{<\lambdach}$) the preimage
$\pi^{-1}(\Gr_{G,\lambdach})$ (resp. $\pi^{-1}(\ol\Gr_{G,\lambdach}),\
\pi^{-1}(\ol\Gr_{G,\lambdach}-\Gr_{G,\lambdach})$).
For $\lambdach\in Y^+$ we denote by $\Lambda_{I,\lambdach}$ (resp. 
$\Lambda_{I,\leq\lambdach},\ \Lambda_{I,<\lambdach}$) the preimage
$\pi_I^{-1}(\Gr_{G,\lambdach})$ (resp. $\pi_I^{-1}(\ol\Gr_{G,\lambdach}),\
\pi_I^{-1}(\ol\Gr_{G,\lambdach}-\Gr_{G,\lambdach})$).
Clearly, $\Lambda_{<\lambdach}$ (resp. $\Lambda_{I,<\lambdach}$)
is closed in $\Lambda_{\leq\lambdach}$ (resp. $\Lambda_{I,\leq\lambdach}$),
with the open complement $\Lambda_\lambdach$ (resp. $\Lambda_{I,\lambdach}$).
In $K$-groups we have exact sequences (see ~\cite{cg} Chapter 5)
$$0\to K^T(\Lambda_{<\lambdach})\to K^T(\Lambda_{\leq\lambdach})\to
K^T(\Lambda_\lambdach)\to0,$$
$$0\to K^T(\Lambda_{I,<\lambdach})\to K^T(\Lambda_{I,\leq\lambdach})\to
K^T(\Lambda_{I,\lambdach})\to0.$$
Thus we obtain a support filtration on $K^T(\Lambda_I)$ (resp.
$K^T(\Lambda)$) with associated graded $\bigoplus_{\lambdach\in Y^+}
K^T(\Lambda_{I,\lambdach})$ (resp. $\bigoplus_{\lambdach\in Y^+}
K^T(\Lambda_\lambdach)$).

We have the embeddings $K^T(\Lambda_\lambdach)\hookrightarrow
K^T(\Lambda_{I,\lambdach})\hookrightarrow K^T(\ul\fu|_{\Gr_\lambdach})\simeq
K^T(\Gr_\lambdach)$. The Weyl group $W$ acts naturally both on 
$K^T(\Lambda_\lambdach)$ and $K^T(\Gr_\lambdach)$, and to prove 
Theorem ~\ref{none} it suffices to check that the image of
$(K^T(\Lambda_\lambdach))^W$ in $K^T(\Lambda_{I,\lambdach})$ contains
(equivalently, coincides with) the intersection 
$K^T(\Lambda_{I,\lambdach})\cap(K^T(\Gr_\lambdach))^W$.

To this end recall that $\Gr_{G,\lambdach}$ can be $G$-equivariantly identified
with the total space $\widetilde\CB$ of a vector bundle 
over a certain partial flag variety
$\CB$ of the group $G$ (the quotient $G/P_\lambdach$ by a parabolic
subgroup depending on $\lambdach$). The Borel subgroup $B\subset G$
acts on $\CB$ with finitely many orbits numbered by the cosets of 
parabolic Weyl subgroup $W^\lambdach=W/W_\lambdach$: we have
$\CB=\bigsqcup_{w\in W^\lambdach}\CB_w$. Let us denote by $\CL\subset T^*\CB$
the union of conormal bundles $\CL=\bigsqcup_{w\in W^\lambdach}T^*_{\CB_w}\CB$.
Let us also denote by $\widetilde\CB_w$ the preimage of $\CB_w$ in 
$\widetilde\CB$ (it coincides with a certain Iwahori orbit 
$\Gr_G^\much\subset\Gr_{G,\lambdach}=\widetilde\CB$).
We define $\widetilde\CL:=\bigsqcup_{w\in W^\lambdach}T^*_{\widetilde\CB_w}
\widetilde\CB\subset T^*\widetilde\CB$.
Then there exists a $G$-equivariant profinite dimensional vector bundle 
$\CV\stackrel{p}{\to}T^*\widetilde\CB$ such that 
$\CV\simeq\ul\fu|_{\Gr_\lambdach}$,
and under this isomorphism we have 
$\CV|_{\widetilde\CL}\simeq\Lambda_{I,\lambdach},\
\CV|_{\widetilde\CB\hookrightarrow T^*\widetilde\CB}\simeq\Lambda_\lambdach$.
Thus to prove Theorem~\ref{none} it is enough to check that the
image of $(K^T(\widetilde\CB))^W$ in $K^T(T^*\widetilde\CB)$ contains
the intersection $K^T(\widetilde\CL)\cap(K^T(T^*\widetilde\CB))^W$.
Equivalently, we have to check that the
image of $(K^T(\CB))^W$ in $K^T(T^*\CB)$ contains
the intersection $K^T(\CL)\cap(K^T(T^*\CB))^W$.
This is the subject of the following lemma.

\begin{lem}
\label{last?}
Let $\imath:\ \CB\hookrightarrow T^*\CB$ denote the embedding of the
zero section, and let $\jmath:\ \CL\hookrightarrow T^*\CB$ denote the
natural closed embedding. Then $\imath_*(K^T(\CB))^W$ coincides with
$\operatorname{Im}
\left(\jmath_*:\ K^T(\CL)\hookrightarrow K^T(T^*\CB)\right)\bigcap
\left(K^T(T^*\CB)\right)^W$.
\end{lem}

\begin{proof}
For $w\in W^\lambdach$ we denote by $w\in\CB_w\subset\CB$ the corresponding
$T$-fixed point. We denote by $i_w$ the closed embedding
$T_w^*\CB\hookrightarrow T^*\CB$
(and also the closed embedding $w\hookrightarrow\CB$, when the confusion is 
unlikely), and we denote by $\fri_w$ the
closed embedding $w\hookrightarrow T^*\CB$. According to Lemma ~\ref{BB},
the image of $\jmath_*:\ K^T(\CL)\hookrightarrow K^T(T^*\CB)$
coincides with the image of $\oplus_{w\in W^\lambdach}i_{w*}:\
\oplus_{w\in W^\lambdach}K^T(T^*_w\CB)\to K^T(T^*\CB)$.
We have an embedding $\oplus_{w\in W^\lambdach}\fri_w^*:\
K^T(T^*\CB)\hookrightarrow\oplus_{w\in W^\lambdach}K^T(w)$, and
similarly an embedding $\oplus_{w\in W^\lambdach}i_w^*:\
K^T(\CB)\hookrightarrow\oplus_{w\in W^\lambdach}K^T(w)$.

Clearly, the $W$-invariants project injectively into any direct summand:
$K^G(\CB)=(K^T(\CB))^W\stackrel{i_w^*}{\hookrightarrow}K^T(w)$
(resp. $K^G(T^*\CB)=(K^T(T^*\CB))^W\stackrel{\fri_w^*}{\hookrightarrow}K^T(w)$)
for any $w\in W^\lambdach$.
Thus it suffices to check that for any $w\in W^\lambdach$ we have a
coincidence
$\operatorname{Im}(\fri_w^*i_{w*}:\ K^T(T^*_w\CB)
^W
\to K^T(w))=
\operatorname{Im}(\fri_w^*\imath_*\on{Res}^G_T:\ K^G(\CB)\to K^T(w))$.
Note that if $w=e$ (the identity coset of $W_\lambdach$ in $W$), then
the image $i_e^*(K^T(\CB))^W\subset K^T(e)$ (resp. 
$\fri_e^*(K^T(T^*\CB))^W\subset K^T(e)$) coincides with 
$(K^T(e))^{W_\lambdach}=\BC[T]^{W_\lambdach}$. Moreover, under identification
$K^T(T^*_e\CB)=K^T(e)=\BC[T]$, we have $K^T(T^*_e\CB)\cap(K^T(T^*\CB))^W=
\BC[T]^{W_\lambdach}$.

Identifying both $K^T(T^*_e\CB)$ and $K^T(e)$ with $\BC[T]$, the map
$\fri_e^*i_{e*}$ is a multiplication by the product
$\Delta_1=\prod_{k=1}^{\dim\CB}(1-\chi_k)$ where $\chi_k$ run through
the characters of $T$ in the tangent space $T_e(T^*_e\CB)=T^*_e\CB$.
Furthermore, identifying $K^G(\CB)$ with $\BC[T]^{W_\lambdach}$, and 
$K^T(e)$ with $\BC[T]$, the map $\fri_e^*\imath_*\on{Res}^G_T$ is a
multiplication by the product
$\Delta_2=\prod_{k=1}^{\dim\CB}(1-\chi'_k)$ where $\chi'_k$ run through
the characters of $T$ in the tangent space $T_e\CB$. We can arrange the
characters $\chi'_k$ so that we have $\chi'_k=\chi^{-1}_k$. Then we see
that $\Delta_1=\Delta_2\cdot\prod_{k=1}^{\dim\CB}(-\chi_k)$, so they
differ by an invertible function, hence the corresponding images coincide:
$\Delta_1\cdot\BC[T]^{W_\lambdach}=\Delta_2\cdot\BC[T]^{W_\lambdach}$.

This completes the proof of the lemma along with Theorem ~\ref{none}.

\subsection{}
\label{classical_langlands}
In this subsection we describe (without striving for high precision)
a conjectural picture motivating Theorem ~\ref{none}.

 We hope that the isomorphism
$K^\GO(\Lambda_G)=\BC[\LT\times T]^W=\BC[T\times\LT]^W=
K^{\LG(\bO)}(\Lambda_\LG)$ lifts to an equivalence of monoidal categories
$F:\ D^bCoh^\GO_{\Lambda_G}(\ul\fu_G)\simeq 
D^bCoh^{\LG(\bO)}_{\Lambda_\LG}(\ul\fu_\LG)$.
The conjectural equivalence $F$ is related to the Langlands correspondence
in the following way.

 Recall that the conjectural (for
$G=GL(n)$
mostly proven   in \cite{denis_pro_langlands}) geometric Langlands 
correspondence is an equivalence of triangulated categories
between the derived category of
$D$-modules on the stack $\on{Bun}_G$ of 
$G$-bundles on a given smooth projective curve $C$,
and the derived category of coherent sheaves on
the stack of $\LG$ local systems on the same curve. One might expect
its ``classical limit'' to be an equivalence between
the derived categories of coherent sheaves  
$L:\ D(T^*\on{Bun}_G)\simeq D(T^*\on{Bun}_\LG)$ where $T^*\on{Bun}_G$
is the cotangent bundle to the moduli stack of $G$-bundles on $C$.
 Given a point $c\in C$, and identifying $\bO$ with the algebra of
functions on the formal neighbourhood of $c$, one gets an action
of $D^bCoh^\GO_{\Lambda_G}(\ul\fu_G)$ on $D(T^*\on{Bun}_G)$. The ``classical
limit'' of 
 the Hecke eigen-property of geometric Langlands correspondence
(see \cite{bd}) should be stated in terms of this action; it should say that
the global equivalence $L$ is compatible with our local equivalence $F$.

\end{proof}

\section{Perverse sheaves and fusion}
\label{seven}

We refer the reader to ~\cite{b} for the definition of perverse equivariant
coherent sheaves and related objects.

\subsection{}
Recall the setup of ~\ref{reminder}. Note that all the
$\GO$-orbits in a connected component of $\Gr_G$ have dimensions of the same
parity. Thus it makes sense to consider the middle perversity function
$p(\Gr_{G,\lambdach})=-\frac{1}{2}
\dim(\Gr_{G,\lambdach})=-\langle\rho,\lambdach\rangle$. It is obviously
strictly monotone and comonotone, but at some connected components of
$\Gr_G$ it takes values in half-integers. This means that we consider
equivariant complexes formally placed in half-integer homological degrees.
The theory of ~\cite{b} defines the artinian abelian category
$\CP^\GO(\Gr_G)$ of perverse $\GO$-equivariant coherent sheaves (with respect
to the above middle perversity).
Let $D^{b,\GO}(\Gr_G)$ denote the bounded derived category of $\GO$-equivariant
coherent sheaves on $\Gr_G$ (with the same convention that the complexes at
``odd'' connected components are placed in half-integer homological degrees).

Given two complexes $\CF,\CG\in D^{b,\GO}(\Gr_G)$ we have their convolution
$\CF\star\CG\in D^{b,\GO}(\Gr_G)$. Recall that
$\CF\star\CG=\Pi_{0*}(\CF\ltimes\CG)$ where $\Pi_0:\ \GF\times_\GO\Gr_G\to\Gr_G$
is the convolution diagram, and $\CF\ltimes\CG$ is the twisted product of
$\CF$ and $\CG$ on $\GF\times_\GO\Gr_G$.

\begin{prop}
\label{preserve}
The convolution preserves perverse sheaves: for $\CF,\CG\in\CP^\GO(\Gr_G)$
we have $\CF\star\CG\in\CP^\GO(\Gr_G)$.
\end{prop}

\begin{proof}
 Denote the projection $\GF\to\GF/\GO=\Gr_G$ by $p$,
and consider a stratification
$\GF\times_\GO\Gr_G=\bigsqcup_{\lambdach,\much\in Y^+}
p^{-1}(\Gr_{G,\lambdach})\times_\GO\Gr_{G,\much}$.
Clearly, $\CF\ltimes\CG$ is smooth (locally free) along this stratification,
and perverse (with respect to the middle perversity).
According to ~\cite{mv} ~2.7, the map $\Pi_0$ is stratified semismall with
respect to the above stratification. Now the perversity of
$\Pi_{0*}(\CF\ltimes\CG)$ follows in the same manner as in the constructible
case, cf. {\em loc. cit.}
\end{proof}

\subsection{The absence of commutativity constraint}
\label{bad luck}
According to Proposition ~\ref{preserve}, $\CP^\GO(\Gr_G)$ acquires the
structure of abelian artinian monoidal category. Moreover, according
to ~\ref{K} ~(a), its $K$-ring is commutative. Nevertheless, $\CP^\GO(\Gr_G)$
admits no commutativity constraint, as can be seen in the following example.

We recall the setup of ~\ref{K-teorija}, and consider $\Gr_{PGL_2}$.
One can check that there are the nonsplit exact sequences in
$\CP^{PGL_2(\bO)}(\Gr_{PGL_2})$:
$$0\to\CV(0)_0\to\CV(0)_1\star\CV(-2)_1\to\CV(-2)_2\to0$$
$$0\to\CV(-2)_2\to\CV(-2)_1\star\CV(0)_1\to\CV(0)_0\to0$$
Thus $\CV(0)_1\star\CV(-2)_1$ and $\CV(-2)_1\star\CV(0)_1$ are nonisomorphic.

\subsection{$\GO\ltimes\Gm$-equivariant sheaves and fusion}
\label{?}
The orbits of $\GO\ltimes\Gm$ on $\Gr_G$ coincide with the $\GO$-orbits,
so one can consider the abelian artinian monoidal category
$\CP^{\GO\ltimes\Gm}(\Gr_G)$ of $\GO\ltimes\Gm$-equivariant coherent perverse
sheaves on $\Gr_G$.
For $\CF\in\CP^{\GO\ltimes\Gm}(\Gr_G)$ we have
$R\Gamma(\Gr_G,\CF)\in D^b(\GO\ltimes\Gm-mod)$.

B.~Feigin and S.~Loktev define (under certain restrictions)
in ~\cite{fl} the {\em fusion product}
$V_1\star\ldots\star V_k\in\GO\ltimes\Gm-mod$ of $\GO\ltimes\Gm$-modules
$V_1,\ldots,V_k$. We recall some of their results in case $G=PGL_2$.

Let $V(n)$ be the $n+1$-dimensional $\GO\ltimes\Gm$-module factoring through
$\GO\ltimes\Gm\twoheadrightarrow G\times\Gm\twoheadrightarrow G$.
Recall the irreducible $PGL_2(\bO)$-equivariant perverse
sheaf $\CV(n)_m$ introduced in ~\ref{K-teorija}. It can be lifted to the
same named $PGL_2(\bO)\ltimes\Gm$-equivariant perverse sheaf, where
the action of $\Gm$ in the fiber over a $\Gm$-fixed point in the orbit
$\Gr_{PGL_2,m}$ is set {\em trivial}.
In particular, $R\Gamma(\Gr_{PGL_2},\CV(n)_1)=V(n)[\frac{1}{2}]$ for $n\geq0$.

Now we can reformulate Theorem ~2.5 of ~\cite{fl} as follows.

\begin{prop}
\label{Feigin}
Let $n_1\geq n_2\geq\ldots\geq n_k$. Then

(a) $R\Gamma(\Gr_{PGL_2},\CV(n_1)_1\star\ldots\star\CV(n_k)_1)$
is concentrated in degree $-\frac{k}{2}$;

(b) $R\Gamma(\Gr_{PGL_2},\CV(n_1)_1\star\ldots\star\CV(n_k)_1)[-\frac{k}{2}]
\simeq V(n_k)\star\ldots\star V(n_1)$.
\end{prop}

\subsection{Multiplication table}
\label{table}
According to Proposition ~\ref{Feigin}, the calculation of fusion product
in $K(\GO\ltimes\Gm-mod)$ is closely related to the ring structure of
$K^{\GO\ltimes\Gm}(\Gr_G)$. Let us formulate the recurrence relations in
$K^{\GO\ltimes\Gm}(\Gr_G)$, compare ~\cite{fl}, end of section ~2.1.
So $\bv(n)_m$ is the class of $\CV(n)_m$ in $K^{\GO\ltimes\Gm}(\Gr_G)$.
We assume that $n\geq0$.

\begin{equation}
\label{odin}
q^{-l}\bv(l+n)_1 0\star\bv(l)_1=q^{-2l}\bv(2l+n)_2+q^2\bv(n-2)_0+q^4\bv(n-4)_0
+\ldots
\end{equation}
(the last summand being $q^n\bv(0)_0$ if $n$ is even,
and $q^{n-1}\bv(1)_0$ if $n$ is odd.)

\begin{equation}
\label{dva}
q^{-l-2}\bv(l-n)_1 0\star\bv(l)_1=
q^{-2l-2}\bv(2l-n)_2+q^{-2}\bv(n-2)_0+q^{-4}\bv(n-4)_0
+\ldots
\end{equation}
(the last summand being $q^{-n}\bv(0)_0$ if $n$ is even,
and $q^{-n+1}\bv(1)_0$ if $n$ is odd.)

\begin{equation}
\label{tri}
\bv(l+1)_1^{\star a}\star\bv(l)_1^{\star b}=
q^{\frac{1}{2}(a(1-a)+l(a+b)(1-a-b))}\bv(a+l(a+b))_{a+b}
\end{equation}

\end{document}